\documentclass[a4paper,11pt]{article} 
\usepackage{amsmath, amssymb, esint, comment}
\usepackage{amscd}
\usepackage{shadow, epsfig, color}
\usepackage{graphicx}
\usepackage{url}
\usepackage{mathrsfs}
\usepackage[final]{changes}
\usepackage{xcolor,cancel}
\usepackage{fancyhdr}

\newcommand{\N}{\mathbb N}

\newcommand{\Q}{\mathbb Q}
\newcommand{\R}{\mathbb R}
\newcommand{\C}{{\sf Ch}}

\renewcommand{\L}{{\rm Lip}}
\newcommand{\e}{\varepsilon}
\newcommand{\1}{\mathbf 1}
\newcommand{\qed}{\ \hfill \fbox{} \bigskip}
\newcommand{\proof}{{\it Proof.} }
\newcommand{\dis}{\displaystyle}

\newcommand{\E}{\mathcal E}

\newcommand{\F}{\mathcal F}
\newcommand{\x}{\overline{x}}

\newcommand{\EN}{\overline{\N}}

\newcommand{\la}{\langle}
\newcommand{\ra}{\rangle}
\newcommand{\m}{\widetilde{m}}

\newcommand{\EE}{\mathbb E}

\renewcommand{\tilde}{\widetilde}
\renewcommand{\b}{\mathbf{b}}

\numberwithin{equation}{section}

\newtheorem{thm}{Theorem}[section]
\newtheorem{defn}[thm]{Definition}
\newtheorem{lem}[thm]{Lemma}
\newtheorem{rem}[thm]{Remark}
\newtheorem{cor}[thm]{Corollary}
\newtheorem{exa}[thm]{Example}
\newtheorem{prop}[thm]{Proposition}
\newtheorem{asmp}[thm]{Assumption}

\makeatletter
\newcommand{\subjclass}[2][2010]{%
  \let\@oldtitle\@title%
  \gdef\@title{\@oldtitle\footnotetext{#1 \emph{Mathematics Subject Classification.} #2}}%
}
\newcommand{\keywords}[1]{%
  \let\@@oldtitle\@title%
  \gdef\@title{\@@oldtitle\footnotetext{\emph{Key Words and Phrases.} #1.}}%
}
\makeatother

\renewcommand{\tilde}{\widetilde}

\title{Convergence of Non-symmetric Diffusion Processes on RCD Spaces 
}
  \author{Kohei Suzuki \thanks{Institute for Applied Mathematics, University of Bonn, Endenicher Allee 60, D-53115 Bonn. E-mail: suzuki@iam.uni-bonn.de}
}

\subjclass{Primary 60F17; Secondary 53C23.}
\keywords{Riemannian Curvature-Dimension Condition, Measured Gromov--Hausdorff Convergence, Non-symmetric Diffusions, Weak Convergence}

\begin{document}
\date{}
\maketitle
\thispagestyle{fancy}
\begin{abstract}
We construct non-symmetric diffusion processes associated with \added{Dirichlet forms consisting of} uniform\added{ly} elliptic forms and derivation operators \added{with killing terms} \replaced{on RCD spaces by aid of non-smooth differential structures {introduced} by Gigli \cite{G16}.}{by aid of non-smooth differential structures (Gigli \cite{G16}) on metric measure spaces under Riemannian Curvature-Dimension conditions.} After constructing diffusions, we investigate conservativeness and the weak convergence of the laws of diffusions in terms of a geometric convergence of the underling spaces and convergences of the corresponding coefficients.
\end{abstract}
\tableofcontents

\section{Introduction}
\subsection{Motivation and Overview}
The aim of this paper is to investigate non-symmetric diffusion processes and their convergence on varying metric measure spaces under Riemannian Curvature-Dimension (RCD) conditions.
We first construct non-symmetric diffusion processes on metric measure spaces under RCD conditions, which are constructed by certain Dirichlet forms consisting of uniformly elliptic operators and derivation operators with killing terms.\
{Then} we investigate conservativeness and the weak convergence of these diffusions in terms of a geometric convergence of the underlying spaces and convergences of the corresponding coefficients.

The notions of CD/RCD conditions on metric measure spaces are generalizations of the notion of lower Ricci curvature bounds in the framework of metric measure spaces, which are stable under geometric convergences such as the measured Gromov--Hausdorff (GH) convergence. They therefore contain various  (finite- and infinite-dimensional) singular spaces such as Ricci limit spaces (Sturm \cite{Sturm06, Sturm06-2}, Lott--Villani \cite{LV09}), Alexandrov spaces (Petrunin \cite{Pet11}, Zhang--Zhu \cite{ZZ10}), warped products and cones (Ketterer \cite{Ket14, Ket14a}), quotient spaces (Galaz-Garc\'ia--Kell--Mondino--Sosa \cite{GKMS17}) and infinite-dimensional spaces such as Hilbert spaces with log-concave measures (Ambrosio--Savar\'e--Zambotti \cite{ASZ09}) (related to various stochastic partial differential equations).
The main point\deleted{ of these notions} is that the notion of lower Ricci curvature bounds can be completely characterized by convexity of \replaced{e}{E}ntropy functionals on Wasserstein spaces, for which only metric measure structures are essential (Sturm \cite{Sturm06, Sturm06-2}, Lott--Villani \cite{LV09}, Ambrosio-Gigli-Savar\'e \cite{AGS14b}, \added{Ambrosio--Gigli--Mondino--Rajala} \cite{AGMR15}, Ambrosio--Mondino-Savar\'e \cite{AMS16} and Erbar--Kuwada--Sturm \cite{EKS15}).

A natural issue in probability theory is whether one can construct diffusion processes on these non-smooth spaces, and if one can construct them, what properties these diffusion processes have. 
By recent developments of analysis on metric measure spaces, we can construct Brownian motions on RCD spaces by using a certain quadratic form, what is called {\it Cheeger energy}. This is a generalization of Dirichlet energy on smooth manifolds and induces a quasi-regular strongly local conservative symmetric Dirichlet form (Ambrosio-Gigli-Savar\'e \cite{AGS14, AGS14b} and \added{Ambrosio--Gigli--Mondino--Rajala} \cite{AGMR15}). Since the Cheeger energy is determined only by the underlying metric measure structure, theoretically speaking, every property of Brownian motions should be derived from the geometric properties of the underlying spaces. With this motivation, in \cite{S17}, the author focused on the relation between geometric and stochastic convergences: the former is  the pointed measured Gromov (pmG) convergence of the underlying spaces, and the latter is the weak convergence of Brownian motions. The main results in \cite{S17} state that the pmG convergence of the underlying spaces implies the weak convergence of Brownian motions on RCD$(K,\infty)$ spaces (and under more strict conditions, these two convergences are equivalent).

In this paper, as a next step of \cite{S17}, we construct non-symmetric diffusion processes and investigate their convergences on varying RCD spaces.  To construct non-symmetric diffusions, we utilize linear transformations between $L^2$-vector fields (called tangent module in Gigli \cite{G16}) as second-order perturbations, and derivation operators (Weaver \cite{W00}) as first-order perturbations corresponding to vector fields on metric measure spaces. We take advantage of a Dirichlet form approach to construct diffusion processes on these non-smooth spaces (see \replaced{Remark \ref{rem: DA} for different approaches}{for the Girsanov transform approach in Fitzsimmons \cite{F07} and the martingale problem approach in Trevisan \cite{T14}}).  
\replaced{Next}{After constructing Dirichlet forms,} we investigate conservativeness and the weak convergence of these diffusion processes. For the weak convergence, we utilize the notion of convergence of non-symmetric forms according to Hino \cite{H98} with a slight modification for varying metric measure spaces. We show the convergence of non-symmetric forms under convergences of uniformly elliptic operators and derivations whereby the convergence of derivations  was introduced by Ambrosio--Stra--Trevisan \cite{AST16}. Consequently, we obtain the weak convergence of the laws of finite-dimensional distributions.  Finally, we study tightness of a sequence of these diffusion processes by aid of the Lyons-Zheng decomposition for non-symmetric forms 
in the case of non-compact spaces. In the case of compact spaces, we use heat kernel estimates. 

We remark that every {result in this paper} can be applied also to the case of time-dependent coefficients $A_t$ (diffusion coefficients), $\b_t$ (drift coefficients) and $c_t$ (killing coefficients) with slight modifications,  but we only deal with the time-independent case in this paper. 



\subsection{Main Results}
\added{In this section, we briefly present our main results, referring to Section 2 for more details on notation.}
\replaced{W}{In this paper,w}e consider {\it a pointed metric measure (p.m.m.)\ space} $\mathcal X=(X, d, m, \x)$, whereby \added{we always assume that}
\begin{align} \label{asmp: basic}
&\text{$(X,d)$ is a complete separable geodesic metric space with non-negative and non-zero Borel} \notag
\\
&\ \text{measure $m$ which is finite on all bounded sets \added{with ${\rm supp}[m]=X$, and $\x$ is a fixed point in $X$}.} 
\end{align}
Here ${\rm supp}[m]$ denotes the support of the measure $m$. 
\added{We also assume the following volume growth condition:
there exist constants $c_1, c_2 >0$ depending only on $K$ satisfying}
\begin{align} \label{asmp: basic1}
m(B_r(\x)) \le c_1e^{c_2r^2}, \quad \forall r>0.
\end{align}
We always assume that $(X,d,m,\x)$ satisfies the RCD$(K,\infty)$ condition, which means that the Ricci curvature is bounded from below by $K$ and the space admits a linear gradient structure in this generality (see Section \ref{subsec: RCDLV}).

We first construct non-symmetric diffusion processes associated with the following bilinear form $\mathcal E: {\rm Lip}_{bs}(X) \times {\rm Lip}_{bs}(X) \to \R$: 
\begin{align} \label{form: non-sym}
\mathcal E (f,g):=\frac{1}{2}\int_X\langle A\nabla f, \nabla g\rangle dm+\int_{X}\b_1(f)gdm+\int_{X}f\b_2(g)dm+\int_Xfgcdm.
\end{align}
In this generality, defining the above formula \eqref{form: non-sym} is a non-trivial issue, but it is possible according to non-smooth differentiable structures developed by Gigli \cite{G16}: ${\rm Lip}_{bs}(X)$ denotes the set of bounded Lipschitz functions with bounded support on $X$, \added{and $|\nabla f|$ denotes the minimal weak upper gradient of $f$}. Let $A: L^{2}(TX) \to L^2(TX)$ denote a (not necessarily symmetric) \replaced{module morphism}{linear operator} on $L^2$-vector field\added{s} (tangent module) so that there exists $H \in L^1_{loc}(X,m)$ satisfying $|A\nabla f| \le H|\nabla f|$ for any $f \in {\rm Lip}_{bs}(X)$. We write $|A|$ for the minimal element among such $H$. Here $|Y| \in L^2(X,m)$ denotes the point-wise norm for $Y \in L^2(TX)$. We denote by $\langle \cdot, \cdot \rangle$ the point-wise scalar product on $L^2(TX)$. The notation $\b_i\ (i=1,2)$ means a derivation operator (see Weaver \cite{W00}, Fitzsimmons \cite{F07}, Gigli \cite{G16}), which is a linear map $\b_i: {\rm Lip}_{bs}(X) \to L^1_{loc}(X,m)$ so that there exists $h \in L^1_{loc}(X,m)$ satisfying $\b_i(f) \le h|\nabla f|$ for any $f \in {\rm Lip}_{bs}(X)$. We write $|\b_i|$ for the minimal element among such $h$. \added{Every notion in this paragraph is explained in Section \ref{sec: Pre} in more detail.} 

Under suitable assumptions (Assumption \ref{asmp: NS1}\deleted{ and \ref{asmp: NS2}}), we show that $\E$ \replaced{is}{becomes} closable and the smallest closed extension $(\E, \F)$ \replaced{is}{becomes} a (quasi-)regular local \deleted{lower-bounded (semi-)}Dirichlet form\deleted{s} (Proposition \ref{thm: NS1}\deleted{, \ref{thm: NS1-2}}). \added{Therefore, there exist diffusion processes corresponding to the Dirichlet form $(\E, \F)$.}

We now focus on the weak convergence of the laws of the corresponding diffusion processes. \added{Let $\mathbb S^\nu=(\mathbb P^{\nu}, S)$ (resp.\ $\hat{\mathbb S}^\nu=(\hat{\mathbb P}^{\nu}, \hat{S})$) denote the diffusion process (resp.\ its dual process) with the initial distribution $\nu$ associated with the Dirichlet form $(\E, \F)$ (resp.\ $(\hat{\E}, \hat{\F})$). Let $\zeta_{\mathbb S}$ and $\zeta_{\hat{\mathbb S}}$ be lifetimes for $\mathbb S$ and $\hat{\mathbb S}$ respectively. Let $\zeta:=\min\{\zeta_{\mathbb S}, \zeta_{\hat{\mathbb S}}\},$ and $\mathbb S_{T}^\nu$ (resp.\ $\hat{\mathbb S}_{T}^\nu$) denotes the diffusion process $\mathbb S^\nu$ $(resp.\ \hat{\mathbb S}^\nu)$ restricted on $\{\zeta>T\}$ for $T>0$.}
We assume the following: 
\begin{asmp} \label{asmp: NS3} \normalfont
Let $\{\mathcal X_n\}_{n \in \N}$ be a sequence of p.m.m.\ spaces satisfying RCD$(K,\infty)$ condition with $m_n(X_n)=1$, or RCD$^*(K,N)$.
Let us suppose the following conditions:
\begin{itemize}
\item[(i)] $\mathcal X_n \to \mathcal X_\infty$ in the pmG sense;
\item[(ii)] $\sup_{n \in \EN}\Bigl(\||A_n|\|_\infty+\||\b^n_1|\|_\infty+\||\b^n_2|\|_\infty+\|{\rm div} \b_1^n\|_\infty+\|{\rm div} \b_2^n\|_\infty+\|c_n\|_\infty\Bigr)<\infty,$ and $A_n$ is symmetric and there exists $\lambda >0$ so that $$\langle A_n \nabla f, \nabla f \rangle \ge \lambda \langle \nabla f, \nabla f \rangle, \ m_n\text{-a.e.},\quad \forall f \in \L_{bs}(X_n), \forall n \in \EN;$$
\item[(iii)] for any non-negative $f \in {\rm Lip}_{bs}(X_n), i=1,2, \forall n \in \EN$,
$$\int_{X_n}(\b^n_i(f)+c_nf)dm_n \ge 0
;$$
\item[(iv)] $A_n \to A_\infty$, and  $\b^n_i \to \b_i^\infty$, ${\rm div}\b_i^n \to {\rm div}\b_i^\infty$ $(i=1,2)$ and $c_n \to c_\infty$ strongly in $L^2$, respectively;
\item[(v)] \replaced{The initial distribution}{Let} $\nu_n \in \mathcal P(X_n)$ \replaced{satisfies}{with} $\nu_n(dx)=\phi_nm_n(dx)$ for $n \in \EN$ {with} $\sup_{n \in \N}\|\phi_n\|_{B_r(\x_n), \infty}<\infty$ and $\phi_n \to \phi_\infty$ weakly in $L^2$.
\end{itemize}
\end{asmp}
The notion of the pmG convergence was introduced by Gigli--Mondino--Savar\'e \cite{GMS13} and is recalled in Section \ref{sec: Pre}. The notion of a convergence of $A_n$ on varying metric measure spaces is introduced in Definition \ref{defn: CONVAN}. The $L^2$-strong convergence of derivation operators $\b_i^n$ was introduced by Ambrosio-Stra-Trevisan \cite{AST16} and the precise definition is recalled in Section \ref{sec: CNSF}. The divergence of a derivation  $\b$ is denoted by ${\rm div} \b$, which is recalled in Section \ref{sec: CNSF} (when we write ${\rm div}\b$, we assume implicitly the existence of ${\rm div}\b$).
We mean $\|\phi_n\|_{B_r(\x_n), \infty}:=\text{ess-}\sup_{x \in B_r(\x_n)}|\phi_n(x)|$, whereby $B_r(\x_n)$ means the open ball centered at $\x_n$ with radius $r$. The notion of $L^p$-convergence of functions on varying metric measure spaces is according to Gigli--Mondino--Savar\'e \cite{GMS13} and stated in Section \ref{sec: CNSF}. The space $\mathcal P(X_n)$ denotes the set of Borel probability measures on $X_n$. \added{Note that, since the ${\rm RCD}(K,\infty)$ condition is stable under the pmG convergence (see \cite[Theorem 7.2]{GMS13}), the limit space $\mathcal X_\infty$ also satisfies the ${\rm RCD}(K,\infty)$ condition. Therefore, the diffusion process associated with $(\E_\infty, \F_\infty)$ and the initial distribution $\nu_\infty$ corresponding to \eqref{form: non-sym} can be defined  also on the limit space $\mathcal X_\infty$ and the corresponding diffusion restricted on $\{\zeta_\infty>T\}$ is denoted by $\mathbb S_{\infty, T}^{\nu_\infty}$ (resp.\ $\hat{\mathbb S}_{\infty, T}^{\nu_\infty}$).}




\deleted{Let $\mathbb S_n^\nu=(\mathbb P_n^{\nu_n}, S^n)$ be a diffusion process with the initial distribution $\nu_n$ associated with the Dirichlet forms $(\E_n, \F_n)$ corresponding to \eqref{form: non-sym} under Assumption \ref{asmp: NS3}.  Let $\zeta_{\mathbb S^n}$ and $\zeta_{\hat{\mathbb S}^n}$ be lifetimes for $\mathbb S^n$ and $\hat{\mathbb S}^n$ respectively. Let $\zeta^n:=\min\{\zeta_{\mathbb S^n}, \zeta_{\hat{\mathbb S}^n}\},$ and $\mathbb S_{n,T}^\nu$ $(resp.\ \hat{\mathbb S}_{n,T}^\nu)$ denotes the diffusion process $\mathbb S_n^\nu$ $(resp.\ \hat{\mathbb S}_n^\nu)$ restricted on $\{\zeta^n>T\}$ for $T>0$.}
\deleted{Note that, since the ${\rm RCD}(K,\infty)$ condition is stable under the pmG convergence} 
\deleted{the limit space $\mathcal X_\infty$ also satisfies the ${\rm RCD}(K,\infty)$ condition. Therefore, diffusion processes associated with $(\E_\infty, \F_\infty)$ corresponding to \eqref{form: non-sym} can be defined  also on the limit space $\mathcal X_\infty$ and the corresponding laws are denoted by $\mathbb S_{\infty,T}^\nu$ $(resp.\ \hat{\mathbb S}_{\infty,T}^\nu)$.}

Under the pmG convergence, we can embed each space $X_n$ to a common ambient space $X$ isometrically and thus, we may consider each $X_n$ to be a subset of $X$. Let $C([0,T];X)$ denote the space of continuous paths from $[0,T]$ to $X$ with uniform topology on compact sets. Now we state the following two main theorems.
\begin{thm} \label{thm: NSC1}
Under Assumption \ref{asmp: NS3}, the laws of $\mathbb S_{n,T}^{\nu_n}$ and $\hat{\mathbb S}_{n,T}^{\nu_n}$ converge weakly to $\mathbb S_{\infty,T}^{\nu_\infty}$ and $\hat{\mathbb S}_{\infty,T}^{\nu_\infty}$, respectively in the space $\mathcal P_{\le 1}(C([0,T]; X))$.
\end{thm}
Here $\mathcal P_{\le 1}(C([0,T];X))$ denotes the set of all Borel sub-probability measures (i.e., measures whose total mass is less than or equal to $1$) on $C([0,\infty);X)$. 
The next theorem requires stronger conditions than Theorem \ref{thm: NSC1}, but the initial distribution can be improved to dirac measures $\delta_{\x_n}$.  
\begin{thm}\label{thm: NSC2} 
Suppose Assumption \ref{asmp: NS3} and RCD$^*(K,N)$ with $\sup_{n \in \N}{\rm diam}(X_n)<\infty$. {If \added{${\rm div}\b^n_1=c_n$ (resp.\ ${\rm div}\b^n_2=c_n$), }\deleted{ and each $n \in \EN$,}}
then the law of\deleted{of} $\hat{\mathbb S}_n^{\x_n}$ (resp.\ ${\mathbb S}_n^{\x_n}$)
converges weakly to  $\hat{\mathbb S}_\infty^{\x_\infty}$ (resp.\ ${\mathbb S}_\infty^{\x_\infty})\deleted{)}$
in $\mathcal P(C([0,\infty); X))$.
\end{thm}

\begin{rem} \normalfont We give two remarks about the main results.
\begin{enumerate}
\item[(i)]
The elliptic constant $\lambda>0$ in (ii) of Assumption \ref{asmp: NS3} needs to be uniform in $n \in \EN$, which is used to prove the convergence of Dirichlet forms and appears in \eqref{ineq: lambdaud} in Section \ref{sec: CNSF}.
\item[(ii)]Under the assumption of Theorem \ref{thm: NSC2}, the underlying space $X_n$ is compact for every $n \in \EN$. 
 \end{enumerate}
\end{rem}
Finally we give a criterion for conservativeness of forms associated with \eqref{form: non-sym} (Proposition \ref{thm: NS2}).

\subsection{Organization of the Paper}
The paper is structured as follows\replaced{.}{:} First, the notation is fixed and preliminary facts are recalled in Section \ref{sec: Pre} (no new results are included): 
basic notations and definitions \added{from metric geometry} (Subsection \ref{subsec: Pre}); \deleted{$L^2$-Wasserstein distance (Subsection \ref{subsec: L2}); }pmG convergence (Subsection \ref{subsec: D}); $L^2$-normed modules (Subsection \ref{subsec: L2T}); Tangent module (Subsection \ref{subsec: Tan}); Dirichlet forms (Subsection \ref{subsec: Der}); RCD$(K,\infty)$ and RCD$^*(K,N)$ spaces (Subsection \ref{subsec: RCDLV}). 
In Section \ref{sec: NSDF}, we prove Proposition \ref{thm: NS1}\deleted{, \ref{thm: NS1-2}} to construct a Dirichlet form corresponding to \eqref{form: non-sym}. 
In Section \ref{sec: CNSF}, we show convergence of non-symmetric forms\deleted{ in terms of convergences of $A_n$ and $\b_i^n$}. 
We first recall $L^p$-convergence of functions on varying metric measure spaces. Secondly, we introduce convergence of $A_n$ and recall a notion of convergence of derivations according to Ambrosio--Stra--Trevisan \cite{AST16}. 
Finally, we show convergence of non-symmetric forms with a modification for varying spaces. 
In Section \ref{sec: CFDDNS}, we prove the weak convergence of finite-dimensional distributions of diffusions under Assumption \ref{asmp: NS3}.
In Section \ref{sec: TGHTNS}, we give proofs for the tightness of diffusions under Assumption \ref{asmp: NS3} and complete the proofs of Theorem \ref{thm: NSC1} and \ref{thm: NSC2}.
In Section \ref{sec: CONSVNS}, we show Proposition \ref{thm: NS2}, which is a criterion for conservativeness. 
Finally in Section \ref{sec: exa}, we give examples for which Assumption \ref{asmp: NS3} is satisfied. 

\section{Notation \& Preliminary Results} \label{sec: Pre}
\subsection{\replaced{Preliminary from Metric Measure Geometry}{Notation}} \label{subsec: Pre}
Let $\N=\{0,1,2,...\}$ and $\overline{\N}:=\N \cup \{\infty\}$ be the set of natural numbers and the set of extended natural numbers, respectively.  
Let $(X,d)$ be a complete separable metric space. We write $B_r(x)=\{y \in X: d(x,y)<r\}$ for an open ball centered at $x \in X$ with radius $r>0$.
By using $\mathscr B(X)$, we denote the family of all Borel sets in $(X,d)$; and by $\mathcal B_b(X)$, the set of real-valued bounded Borel-measurable functions on $X$.
Let $C(X)$ denote the set of real-valued continuous functions on $X$, while $C_b(X), C_0(X)$ and $C_{bs}(X)$ denote the subsets of $C(X)$ consisting of bounded functions, functions with compact support, and bounded functions with bounded support, respectively. \added{Let ${\rm Lip}(X)$ denote the set of real-valued Lipschitz continuous functions on $X$. Let ${\rm Lip}_b(X)$ and $\mathrm{Lip}_{bs}(X)$ denote the subsets of ${\rm Lip}(X)$ consisting of bounded functions, and bounded functions with bounded supports, respectively. For $f \in {\rm Lip}_{bs}(X)$, {\it the global Lipschitz constant ${\rm Lip}(f)$} is defined as the infimum of $L>0$ satisfying $|f(x)-f(y)| \le Ld(x,y)$ for any $x,y \in X$.}
The set $\mathcal P(X)$ denotes all Borel probability measures on $X$.
The set of continuous functions on $[0,\infty)$ valued in $X$ is denoted by $C([0,\infty); X)$. 

Let $\mathrm{supp}[m]=\{x \in X: m(B_r(x))>0, \ \forall r>0\}$ denote the support of $m$.
Let $(Y,d_Y)$ be a complete separable metric space. For a Borel measurable map $f: X \to Y$, let $f_\#m$ denote the push-forward measure on $Y$:
$${f}_\#m(B) = m(f^{-1}(B)) \quad \text{for any Borel set} \quad  B \in \mathscr B(Y).$$
\added{For a measure space $(X, m)$ with a Borel measure $m$, we denote by $L^p(X,m)$ ($L^p(m)$ for brevity if no confusion occurs) ($1 \le p \le \infty$) the space of $m$-equivalence classes of Borel measurable functions $f:X \to \R\cup\{\infty\}$ so that $\|f\|_{L^p(X,m)}^p:=\int_X|f|^pdm<\infty$ if $1 \le p <\infty$, and $\|f\|_{L^\infty(X,m)}={\text {ess-sup}}_{x \in X}|f(x)|<\infty$ in the case of $p=\infty$. We sometimes write $\|\cdot\|_p$ for brevity. 
Let $L^0(X,m)$ denote the set of equivalent classes of Borel measurable functions $f: X \to \R$. 
For $f,g \in L^2(X,m)$, let $(f,g)_{L^2(X,m)}$ (simply $(f,g)$) denote the inner product $\int_Xfg dm.$ For a measurable set $A\subset X$, let us denote the indicator function by $\1_A$, which is equal to $1$ for $x \in A$ and $0$ otherwise. For any two functions $f, g: X \to \R$, we write $f\vee g=\max\{f,g\}$ and $f \wedge g=\min\{f,g\}.$}
\deleted{A continuous curve $\gamma: [a,b] \to X$ is said to be {\it connecting $x$ and $y$} if $\gamma_a=x$ and $\gamma_b=y$.
A continuous curve $\gamma: [a,b] \to X$ is said to be {\it a minimal geodesic} if }
\deleted{In particular, if $\frac{d(\gamma_a, \gamma_b)}{|b-a|}$ can be replaced by $1$, we say that $\gamma$ is {\it  unit-speed}.}

A curve $\gamma:[0,1] \to X$ is {\it absolutely continuous} if there exists a function $f \in L^1(0,1)$ so that 
\begin{align} \label{ineq: ABC}
d(\gamma_t, \gamma_s) \le \int_t^sf(r)dr, \quad \forall t,s \in [0,1], \quad t<s.
\end{align}
The metric speed $t \mapsto |\dot{\gamma}|_t \in L^1(0,1)$ is defined as the essential infimum among all the functions $f$
satisfying \eqref{ineq: ABC}.

A Borel probability measure $\pi$ on $C([0,1];X)$ is {\it a test plan} if there exists a constant $C(\pi)$ so that 
$$(e_t)_{\#}\pi \le C(\pi)m, \quad \forall t \in [0,1], \quad \text{with} \quad \int\int_0^1|\dot{\gamma}_t|^2dt d\pi(\gamma)<\infty.$$
Here $e_t(\gamma):=\gamma(t) \in X$ is the evaluation map.

{\it The set of Sobolev functions} $S^2(X,d,m)$ (or, simply $S^2(X)$) is defined to be the space of all functions in $L^0(X,m)$ so that there exists a non-negative $G \in L^2(m)$ for which it holds
$$\int|f(\gamma_1)-f(\gamma_0)|d\pi(\gamma) \le \int\int_0^1G(\gamma_t)|\dot{\gamma}_t|dt d\pi(\gamma), \quad \forall \text{test plan } \pi.$$ 
It turns out (see \cite{AGS14b}), that for $f \in S^2(X)$ there exists a minimal $G$ in the $m$-a.e.\ sense for which the above inequality holds. We denote by $|\nabla f|$ such $G$ and call it {\it minimal weak upper gradient}. 
Let us define $W^{1,2}(X,d,m):=S^2(X,d,m) \cap L^2(X,m)$ (or, simply $W^{1,2}(X)$). A functional {\it Cheeger energy} ${\sf Ch}: W^{1,2}(X,d,m) \to \R$ is defined as follows
\begin{align*}
{\sf Ch}(f)&=\frac{1}{2}\int_X|\nabla f|^2dm, \quad f \in W^{1,2}(X,d,m).
\end{align*}
\begin{rem}\normalfont \label{rem: IH}
Note that ${\sf Ch}: L^2(X,m) \to [0,+\infty]$ is a lower semi-continuous and convex functional, but not necessarily a quadratic form.
This means that $(W^{1,2}(X,d,m), \sqrt{2{\sf Ch}(\cdot)+\|\cdot\|^2_2})$ is a Banach space, but not necessarily a Hilbert space. 
\end{rem}
We say that $(X,d,m)$ satisfies {\it the infinitesimal Hilbertian (IH) condition} if ${\sf Ch}$ is a quadratic form, i.e.,
\begin{align}
	2{\sf Ch}(u)+2{\sf Ch}(v)={\sf Ch}(u+v) + {\sf Ch}(u-v), \label{defn: Cheeger}
	\end{align}
	for any $u,v \in W^{1,2}(X,d,m).$
Let us define the point-wise scalar product as follows
\begin{align} \label{def: PWS}
\langle \nabla f, \nabla g \rangle:=\lim_{\e\to 0}\frac{|\nabla(f+\e g)|^2-|\nabla f|^2}{\e}, \quad f,g \in W^{1,2}(X,d,m),
\end{align}
whereby  the limit is with respect to $L^1(m)$. If the Cheeger energy ${\sf Ch}$ is quadratic, the point-wise inner product becomes a $L^1(m)$-valued bilinear form (see \cite[Definition 4.12]{AGS14b}, and \cite[Theorem 2.7]{AST16}). We set ${\sf Ch}(f,g):=(1/2)\int_X \langle \nabla f, \nabla g \rangle dm.$

\subsection{Pointed Measured Gromov Convergence} \label{subsec: D}
We recall the definition of pmG convergence introduced in Gigli-Mondino-Savar\'e \cite{GMS13}.
\begin{defn}[\cite{GMS13}]\normalfont ({\bf pmG Convergence}) \label{prop: Dconv} 
A sequence of p.m.m.\ spaces $\mathcal X_n=(X_n, d_n, m_n, \x_n)$ satisfying \eqref{asmp: basic} is said to be {\it convergent to $\mathcal X_\infty=(X_\infty, d_\infty, m_\infty, \x_\infty)$ in the pointed measured Gromov (pmG) sense} if there exist a complete separable metric space $(X,d)$ and isometric embeddings $\iota_n: X_n \to X\  (n \in \EN)$
satisfying
\begin{align}\label{eq: VGC}
\iota_n(\x_n) \to \iota_\infty(\x_\infty) \in {X_\infty}, \quad \text{and}\quad \int_{X}f \ d({\iota_n}_\#m_n) \to \int_{X}f \ d({\iota_\infty}_\#m_\infty),
\end{align}
for any bounded continuous function $f:X \to \R$ with bounded support. 
\end{defn}
\begin{rem} \normalfont \label{rem: GHMG}
We give two remarks for Definition \ref{prop: Dconv}.
\begin{enumerate}
\item[(i)] The pmG convergence is weaker than the pointed measured Gromov-Hausdorff (pmGH) convergence (\cite[\replaced{Theorem}{Proposition} 3.30, Example 3.31]{GMS13}). If \deleted{${\rm supp}[m_\infty]=X_\infty$ and }$\{\mathcal X_n\}_{n \in \N}$ satisfies a uniform doubling condition, then pmG and pmGH coincide \cite[\replaced{Theorem}{Proposition} 3.33]{GMS13}.
\item[(ii)] The  pmG convergence is metrizable by a distance $p\mathbb G_W$ on the collection $\mathbb X$ of all isomorphism classes of p.m.m.\ spaces (\cite[Definition 3.13]{GMS13}). The space $(\mathbb X, p\mathbb G_W)$ \replaced{is}{becomes} a complete and separable metric space (\cite[Theorem 3.17]{GMS13}). 
\end{enumerate}
\end{rem}

\subsection{$L^p(m)$-normed Module} \label{subsec: L2T}
In this subsection, we recall the notion of $L^p$-normed module by following \cite[\S 1.2]{G16}. 
\begin{defn}[\cite{G16}] \normalfont {\bf ($L^\infty(m)$-premodule)}  
{\it An $L^\infty(m)$-premodule} is a Banach space $(\mathcal M, \|\cdot\|_\mathcal M)$ equipped with a bilinear map
$L^\infty(m) \times \mathcal M \ni (f,v) \mapsto f\cdot v \in \mathcal M$ satisfying 
$$(fg)\cdot v = f \cdot(g\cdot v), \quad \1\cdot v = v, \quad \|f\cdot v\|_\mathcal M \le \|f\|_{L^\infty(m)}\|v\|_\mathcal M,$$
for any $v \in \mathcal M$ and $f,g \in L^\infty(m)$, whereby $\1:=\1_X \in L^\infty(m)$.
\end{defn}
\begin{defn}[\cite{G16}] \normalfont {\bf ($L^\infty(m)$-module/Hilbert Module)}  
{\it An $L^\infty(m)$-module} is an $L^\infty(m)$-premodule $\mathcal M$ satisfying the following two conditions:
\begin{description}
\item[(i)] ({\bf Locality}) For any $v \in \mathcal M$, $A_n \in \mathscr B(X)$ and $n \in \N$, 
$$\1_{A_n}\cdot v=0, \quad \forall n \in \N \quad \text{implies} \quad \1_{\cup_{n \in \N}A_n}\cdot v=0.$$
\item[(ii)] ({\bf Gluing}) For any sequence $\{v_n\}_{n \in \N} \subset \mathcal M$ and $\{A_n\}_{n \in \N} \subset \mathscr B(X)$ so that 
$$\1_{A_i \cap A_j}\cdot v_i=\1_{A_i \cap A_j}\cdot v_j, \quad \forall i,j \in \N, \quad \text{and} \quad \limsup_{n \to \infty} \|\sum_{i=1}^n\1_{A_i}\cdot v_i\|_\mathcal M <\infty,$$
there exists $v \in \mathcal M$ so that 
$$\1_{A_i}\cdot v=\1_{A_i}\cdot v_i, \quad \forall i \in \N, \quad \text{and} \quad \|v\|_\mathcal M \le \liminf_{n \to \infty}\|\sum_{i=1}^n\1_{A_i}\cdot v_i\|_\mathcal M.$$
\end{description}
If furthermore $(\mathcal M, \|\cdot\|_{\mathcal M})$ is a Hilbert space, then $\mathcal M$ is called {\it Hilbert module}.
\end{defn}
\begin{exa} \normalfont
One of the typical examples for $L^\infty(m)$-modules is $\mathcal M=L^p(X, m)$ with the norm $\|\cdot\|_\mathcal M:=\|\cdot\|_{L^p(m)}$ for $p \in [1,\infty]$. If the underlying space $(X,d,m)$ is a Riemannian manifold, then an $L^p(X,m)$-vector field for $p \in [1,\infty]$ is also an $L^\infty(m)$-module. 
\end{exa}

For two given $L^\infty(m)$-modules, $\mathcal M_1, \mathcal M_2$, a map $T: \mathcal M_1 \to \mathcal M_2$ is called {\it a module morphism} provided that it is a bounded linear map from $\mathcal M_1$ to $\mathcal M_2$ as a map between Banach spaces and satisfies 
\begin{align} \label{prop: mod}
T(f \cdot v)=f\cdot T(v), \quad \forall v \in \mathcal M_1, f \in L^\infty(m).
\end{align}
The set of all module morphisms is denoted by $\text{\sc Hom}(\mathcal M_1, \mathcal M_2)$. It is known that $\text{\sc Hom}(\mathcal M_1, \mathcal M_2)$ has a canonical $L^\infty(m)$-module structure. 
\begin{defn}[\cite{G16}] \normalfont {\bf (Dual Module)} 
For an $L^\infty(m)$-module $\mathcal M$, {\it the dual module $\mathcal M^*$} is defined as $\text{\sc Hom}(\mathcal M, L^1(m))$.
\end{defn}
\begin{defn}[\cite{G16}] \normalfont {\bf ($L^p(m)$-normed module)} \label{defn: Lpnormed}
 Let $p \in [1,\infty]$.
An $L^p(m)$-normed module is an $L^\infty(m)$-module $\mathcal M$ endowed with a map $|\cdot|: \mathcal M \to L^p(m)$ with non-negative values so that 
$$\||v|\|_{L^p(m)}=\|v\|_\mathcal M, \quad |f\cdot v|=|f||v|, \quad \text{$m$-a.e,}$$
for every $v \in \mathcal M$ and $f \in L^\infty(m)$. The map $|\cdot|$ is called {\it point-wise norm}.
\end{defn}

\subsection{Tangent Module} \label{subsec: Tan}
In this subsection, following \cite[\S 2]{G16}, we recall the tangent module $L^2(TX)$ on $(X,d,m)$, which is an $L^2(m)$-normed module in the sense of Definition \ref{defn: Lpnormed} and  a generalized notion of the space of $L^2$-sections of the tangent bundle on smooth manifolds. 

We recall the set {\it Pre-cotangent module} Pcm, which is defined as follows:
\begin{defn}[\cite{G16}] \normalfont {\bf (Pre-cotangent module)}
The set {\rm Pcm} defined as follows is called {\it Pre-cotangent module}:
\begin{align*}
{\rm Pcm}:=\Bigl\{\{(f_i, A_i)\}_{i \in \N}: &\{A_i\}_{i \in \N} \subset \mathscr B(X)\ \text{is a partition of X}
\\
&f_i \in S^2(X), \forall i \in \N,\ \text{{and}} \ \sum_{i \in \N}\int_{A_i}|\nabla f|^2dm<\infty \Bigr\}.
\end{align*}
An equivalence relation between two elements in Pcm $\{(f_i, A_i)\}_{i \in \N} \sim \{(g_j, B_j)\}_{j \in \N}$ is defined as follows:
$$|\nabla(f_i-g_j)|=0, \quad m\text{-a.e. on}\ A_i \cap B_j, \quad \forall i,j \in \N.$$
\end{defn}

A vector space structure can be endowed with the quotient space Pcm$/\sim$ by defining the sum and the scalar multiplication as follows: for any $\lambda \in \R$, 
$$[(f_i, A_i)_i]+[(g_j, B_j)_j]:=[(f_i+g_j, A_i\cap B_j)_{i,j}], \quad \lambda[(f_i, A_i)_i]:=[(\lambda f_i, A_i)_i].$$

The product operation $\cdot: {\rm Sf}(m) \times {\rm Pcm}/\sim \to {\rm Pcm}/\sim$ can be defined by the following manner. Let ${\rm Sf}(m) \subset L^\infty(m)$ denote the set of all simple functions $f$, which means that $f$ attains only a finite set of values. Given $[(f_i, A_i)_i] \in {\rm Pmc}/\sim$ and $h=\sum_{j}a_j\1_{B_j} \in {\rm Sf}(m)$ with $\{B_j\}_{j \in \N}$ a partition of $X$, the product $h \cdot [(f_i, A_i)_i]$ is defined as follows:
$$h \cdot [(f_i, A_i)_i]:=[(a_jf_i, A_i\cap B_j)_{i,j}].$$

We now recall the point-wise norm $|\cdot|_*$ (we use the notation $|\cdot|_*$ as a point-wise norm for the sake of consistency with the definition of tangent modules given later): Define $|\cdot|_*$ on ${\rm Pcm}/\sim \to L^2(X,m)$ by
$$\bigl| [(f_i, A_i)_i]\bigr|_*:=|\nabla f_i|, \quad \text{$m$-a.e. on $A_i$ for all $i \in \N.$}$$
Then the map $\|\cdot\|_{L^2(T^*X)}: {\rm Pcm}/\sim \to [0,\infty)$ is defined as follows:
$$\|[(f_i, A_i)_i]\|_{L^2(T^*X)}^2:=\int\bigl| [(f_i, A_i)_i]\bigr|^2_*dm=\sum_{i \in \N}\int_{A_i}|\nabla f_i|^2dm.$$
Then $\|\cdot\|_{L^2(T^*X)}$ is a norm on Pcm$/\sim.$
\begin{defn}[\cite{G16}] \normalfont {\bf (Cotangent Module)} \label{defn: cotan}
The cotangent module $(L^2(T^*X), \|\cdot\|_{L^2(T^*X)})$ is defined as the completion of $({\rm Pcm}/\sim, \|\cdot\|_{L^2(T^*X)})$.
\end{defn}
It can be checked that the cotangent module $(L^2(T^*X), \|\cdot\|_{L^2(T^*X)})$ is an $L^2$-normed module with the product $\cdot$ (which can be extended to the map $\cdot: L^\infty(m) \times L^2(T^*X) \to L^2(T^*X)$), and the point-wise norm $|\cdot|_*$ (see \cite[\S 2.2]{G16} for more details).
\begin{defn}[\cite{G16}] \normalfont {\bf (Tangent Module)} \label{defn: tan}
The tangent module $(L^2(TX), \|\cdot\|_{L^2(TX)})$ is defined as the dual module of $(L^2(T^*X), \|\cdot\|_{L^2(T^*X)})$.
The point-wise norm associated with the dual of $|\cdot|_*$ is written as $|\cdot|$.
\end{defn}
Under the condition (IH), the tangent module $(L^2(TX), \|\cdot\|_{L^2(TX)})$ is a Hilbert module and the point-wise norm $|\cdot|$ satisfies the parallelogram identity. Therefore, we can define the point-wise inner product $\langle \cdot, \cdot \rangle$.

Now we recall the notions of differential and gradient for a function in Sobolev class.
\begin{defn} [\cite{G16}] \normalfont {\bf (Differential)} \label{defn: dif}
Let $f \in S^2(X)$. The differential $df \in L^2(T^*X)$ is defined as 
$$df:=[(f,X)] \in {\rm Pcm}/\sim \subset L^2(T^*X).$$
Here $[(f,X)]$ means $[(f_i,A_i)_{i \in \N}]$ for $f_0=f$, $A_0=X$ and $A_i=\emptyset$ for $i \ge 1$.
\end{defn}
By definition, we have $|df|_*=|\nabla f|$. The notion of gradient of a Sobolev function is defined through duality with the notion of the differential. 
\begin{defn} [\cite{G16}] \normalfont {\bf (Gradient)} \label{defn: grad}
Let $f \in S^2(X)$. We say that $X \in L^2(TX)$ is {\it a gradient {of $f$}} if 
$$df(X)=|X|^2=|df|^2_*.$$
The set of all gradients of $f$ is denoted by ${\rm Grad}(f)$.
\end{defn}
Under condition (IH), the set ${\rm Grad}(f)$ has a unique element, which is denoted by $\nabla f$. In this case, the gradient $\nabla f$ satisfies the following linearity (\cite[Proposition 2.3.17]{G16}):
$$\nabla(f+g)=\nabla f+\nabla g, \quad m\text{-a.e.}, \quad f,g \in S^2(X).$$
Let $(X,d,m)$ satisfy (IH) and  $\langle \cdot, \cdot \rangle: L^2(TX) \to L^1(X,m)$ be the point-wise inner product, which is induced by the structure of $L^2$-normed module with the point-wise norm $|\cdot|$ in $L^2(TX)$. Under condition (IH), $\langle \nabla f, \nabla g \rangle$ can be identified in the $m$-a.e.\ sense with the same expression defined in \eqref{def: PWS} in Subsection \ref{subsec: Pre}.

\subsection{Derivation} \label{subsec: Der}
In this subsection, we briefly explain derivations on metric measure spaces by following \cite[\S 3]{AST16}.
\begin{defn}[\cite{AST16}] {\bf (Derivation)} \normalfont
A linear functional $\b: {\rm Lip}_{bs}(X) \to L^0(m)$ is said to be {\it a derivation} if there exists $h \in L^0(m)$
 so that 
 $$|\b(f)| \le h|\nabla f|,\quad m\text{-a.e. in $X$, for all $f \in {\rm Lip}_{bs}(X).$}$$
 The $m$-a.e.\ smallest function $h$ satisfying the above inequality is denoted by $|\b|$. 
 The space of all derivations is denoted by ${\rm Der}(X,d,m)$. We denote by ${\rm Der}^p(X,d,m)$ (resp.\ ${\rm Der}^p_{loc}(X,d,m)$) the space of derivations $\b$ so that $|\b| \in L^p(m)$ (resp.\ $L^p_{loc}(m)$).
 \end{defn}
Derivation operators satisfy the local property: for any $f,g \in {\rm Lip}_{bs}(X)$,
$$|\b(f-g)| \le h|\nabla (f-g)|=0, \quad m\text{-a.e. on }\ \{f=g\}.$$
By the local property, the chain rule holds
$$\b(\phi(f))=(\phi' \circ f)\b(f), \quad \phi \in {\rm Lip}(\R), \quad m\text{-a.e.},$$
and the Leibniz rule also holds:
$$\b(fg)=\b(f)g+f\b(g), \quad m\text{-a.e.}.$$
See \cite[\S 2.2, 2.3]{G16} and \cite[\S 3]{AST16} for more details.
\begin{rem}\normalfont {\bf (${\rm Der}^2(X,d,m)$ and $L^2(TX)$)}
The space ${\rm Der}^2(X,d,m)$ is an $L^2(m)$-normed module with the map $\b \to |\b|.$
If $W^{1,2}(X,d,m)$ is reflexive as a Banach space, then ${\rm Der}^2(X,d,m)$ can be canonically and isometrically identified with the tangent module $L^2(TX)$ recalled in Definition \ref{defn: tan}. See \cite[Remark 3.5]{AST16} for more details. 
\end{rem}

Now we recall the notion of {\it divergence of derivations}. 
\begin{defn}[\cite{AST16}] {\bf (Divergence)} \normalfont
A derivation $\b \in {\rm Der}^1_{loc}(X,d,m)$ has divergence in $L^1_{loc}(X, m)$ if there exists $g \in L^1_{loc}(X,m)$ so that 
$$-\int_{X}\mathbf b(f)dm=\int_X fg dm, \quad \forall f \in {\rm Lip}_{bs}(X).$$  
Such a $g$ is uniquely determined if it exists, and we denote it by ${\rm div}\mathbf b.$ The existence of such $g$ is not necessarily true for general $\b$ but when we write ${\rm div}\b$, we implicitly assume the existence of such $g$. Let ${\rm Div}^p_{loc}(X,d,m):=\{\b \in {\rm Der}^p_{loc}(X,d,m): {\rm div} \b \in L^p_{loc}(X,d,m)\}$ and ${\rm Div}^p(X,d,m):=\{\b \in {\rm Der}^p_{loc}(X,d,m): {\rm div} \b \in L^p(X,d,m)\}$.
\end{defn}
By using the Leibniz rule, we have 
\begin{align} \label{eq: DIV}
\int_X\mathbf b(f)\phi dm=-\int_X \mathbf b(\phi)fdm-\int_{X}f\phi {\rm div}\mathbf b dm,\quad \forall f, \phi \in {\rm Lip}_{bs}(X).
\end{align}

\subsection{Dirichlet Forms} \label{subsec: NSD}
In this subsection, following \cite{MR92}, we recall basic notions concerning Dirichlet forms.

Let $\mathcal F \subset L^2(X,m)$ be a dense linear subspace and $\mathcal E$ be a bilinear form on $\mathcal F$. We write $\E_{\alpha}(f,g):=\E(f,g)+\alpha(f,g)_{L^2(X,m)}$ and $\E_{\alpha}(f):=\E_\alpha(f,f)=\E(f,f)+\alpha\|f\|_2^2$ for $\alpha \in [0,\infty)$ in short. {\it The symmetric part of $\E$} is defined by $\tilde{\E}(f,g)=(1/2)(\E(f,g)+\E(g,f))$ and {\it the anti-symmetric part of $\E$} by $\check{\E}(f,g)=(1/2)(\E(f,g)-\E(g,f))$. The bilinear form $(\E, \F)$ is {\it a coercive closed form} if the following three conditions hold:
\begin{description}
\item[($\mathcal E.1$)] {\bf $\mathcal E$ is non-negatively definite}: $\E(f) \ge 0$ for all $f \in \mathcal F$.
\item[($\mathcal E.2$)] {\bf $\mathcal E$ satisfies the weak sector condition}: there exists a constant $C \ge 1$ so that 
$$|\E_1(f,g)| \le C \E_{1}(f)\E_{1}(g), \quad \forall f, g \in \mathcal F.$$
\item[($\mathcal E.3$)] $\F$ is a Hilbert space with respect to the symmetric part $\tilde{\E}^{1/2}_1.$
\end{description}
Let $D \subset L^2(X,m)$ be a dense linear subspace. A bilinear form $(\E, D)$ satisfying $(\E.1)$ and $(\E.2)$ is {\it closable} if, for any $\E$-Cauchy sequence $f_n \in D$ with $\lim_{n \to \infty}\|f_n\|_{L^2(m)}=0$, it holds that $\lim_{n \to \infty}\E(f_n)=0.$
We say that $(\E, \F)$ is {\it symmetric} if $\E(f,g)=\E(g,f)$ for all $f,g\in \mathcal F$. {\it The dual form} $\hat{\E}$ is defined to be $\hat{\E}(f,g)=\E(g,f)$ for $f, g \in \mathcal F$. 

If $(\E, \F)$ is a coercive closed form, then there exist the corresponding semigroups $\{T_t\}_{t \ge 0}$ and $\{\hat{T}_t\}_{t \ge 0}$ on $L^2(X,m)$ so that $(T_tf, g)=(g, \hat{T}_tf)$ for any $t \ge 0$ and $f, g \in L^2(X,m)$, and the corresponding resolvents $G_\alpha$ and $\hat{G}_\alpha$, which are defined as $G_\alpha f=\int_0^\infty e^{-\alpha t}T_tf dt$ and $\hat{G}_\alpha f=\int_0^\infty e^{-\alpha t}\hat{T}_t f dt$, satisfy
$$\E_\alpha (G_\alpha f, g)=(f,g)=\E_\alpha(g, \hat{G}_\alpha f), \quad \forall f \in L^2(X,m),\ g \in \mathcal F,\ \alpha>0.$$

Concerning the Markovian property, the following statements are known to be equivalent (e.g., \cite[Theorem 1.1.5.]{O13}):
\begin{description}
\item[($\mathcal E.4$)] For all $f \in \mathcal F$, it holds that 
\begin{align} \label{Markovian}
f^+\wedge 1 \in \mathcal F, \quad \E(u+u^+\wedge1, u-u^+\wedge1) \ge 0, \quad {\E}(u-u^+\wedge1, u+u^+\wedge1) \ge 0.
\end{align}
Here $u^+:=u\vee 0.$
\item[(M)] $\{T_t\}_{t \ge 0}$ and $\{\hat{T}_t\}_{t \ge 0}$ are Markovian: If $f \in L^2(X,m)$ satisfies $0 \le f \le 1$ $m$-a.e., then $0 \le T_tf \le 1$ and  $0 \le \hat{T}_tf \le 1$ $m$-a.e.. 
\end{description}
A bilinear form $(\E, \F)$ is called {\it Dirichlet form} if ($\mathcal E.1$)--($\mathcal E.4)$ hold. 

Now we recall the property of regularity/quasi-regularity for Dirichlet forms, which is a sufficient condition for the existence of Hunt processes/$m$-tight special standard processes and their dual processes (see \cite[Theorem 3.5 Chapter IV]{MR92}) corresponding to Dirichlet forms. An increasing sequence $\{E_n\}_{n \in \N}$ of closed subsets of $X$ is called {\it$\E$-nest} if 
$$\bigcup_{n \in \N }\mathcal F|_{E_n} \text{ is dense in $\mathcal F$ with respect to $\tilde{\E}_1^{1/2}$}.$$ 
Here we mean that $\mathcal F|_{A}:=\{u \in \mathcal F: u=0\ m\text{-a.e.\ on $A^c$}\}$.
 A subset $N \subset X$ is called {\it $\E$-exceptional} if 
 $$N \subset \bigcap_{n \in \N}E_n^c \text{ for some $\E$-nest $\{E_n\}_{n \in \N}$.}$$
  We say that a property of points in $X$ holds $\E$-quasi-everywhere ($\E$-q.e.) if the property
 holds outside some $\E$-exceptional set. A function $f$ $\E$-q.e.\ defined on $X$ is called {\it $\E$-quasi-continuous} if there exists an $\E$-nest $\{E_n\}_{n \in \N}$ so that $f \in C(\{E_n\})$ whereby 
 $$C(\{E_n\}):=\{f: A \to \R: \bigcup_{n \in \N} E_n\ \subset A \subset X,\ f|_{E_n}\text{ is continous} \quad \forall n \in \N\}.$$
A Dirichlet form $(\E, \F)$ on $L^2(X,m)$ is called {\it quasi-regular} if the following three conditions hold:
\begin{enumerate}
\item[(i)] There exists an $\E$-nest $\{E_n\}_{n \in \N}$ consisting of compact sets.
\item[(ii)] There exists an $\tilde{\E}_1^{1/2}$-dense subset of $\F$ whose elements have $\E$-quasi-continuous $m$-versions.
\item[(iii)] There exist $u_n \in \F$ for $n \in \N$ having $\E$-quasi-continuous $m$-versions $\tilde{u}_n$ and an $\E$-exceptional set $N \subset X$ so that $\{\tilde{u}_n\}_{n \in \N}$ separates points of $X \setminus N$.
\end{enumerate}
Let $(X,d)$ be a locally compact separable metric space with a Radon measure $m$. A Dirichlet form $(\E, \F)$ on $L^2(X,m)$ is called {\it regular} with a core $\mathcal C_1$ if $\mathcal C_1 \subset C_0(X)\cap \mathcal F$ is dense both in $C_0(X)$ with the uniform norm $\|\cdot\|_\infty$ and in $\mathcal F$ with $\tilde{\E}^{1/2}_{1}$, respectively. We note that $(\E, \F)$ is quasi-regular if it is regular (\cite[Chapter IV Section 4 a)]{MR92}.

Let $\{T_{t}\}_{t \ge 0}$ be the semigroup corresponding to $(\E, \F)$. An important fact (\cite[Theorem 3.5 Chapter IV]{MR92}) is that if a Dirichlet form $(\E, \F)$ is quasi-regular, then there exists an $m$-tight special standard process $(\Omega, \mathcal M, \{\mathcal M_t\}_{t \ge 0}, \{S_t\}_{t \ge 0}, \{\mathbb P^x\}_{x \in X})$ (\cite[Definition 1.13 in Chapter IV]{MR92}) so that, for all $t \ge 0$ and $f \in \mathcal B_b(X) \cap L^2(X,m)$,
$$T_t f(x)=\mathbb E^x(f(S_t)), \quad \text{$\E$-q.e.\ $x$}.$$
Here $\mathbb E^x(f(S_t)):=\int_\Omega f(S_t(\omega))\mathbb P^x(d\omega)$.

We adjoin an extra point $\partial$ (the cemetery point) to $X$ as an isolated point to obtain a Hausdorff topological space $X_\partial$ with Borel $\sigma$-algebra $\mathscr B(X_\partial)=\mathscr B(X) \cup \{B \cup \{\partial\}: B \in \mathscr B(X)\}.$ Any function $f: X \to \R$ can be considered as a function from $X_\partial$ by defining $f(\partial)=0$. If $X$ is locally compact, we consider the one-point compactification (Alexandroff compactification) for $X_\partial$. We say that a stochastic process $(\Omega, \mathcal M, \{\mathcal M_t\}_{t \ge 0}, \{S_t\}_{t \ge 0}, \{\mathbb P^x\}_{x \in X})$ {\it has a lifetime $\zeta$} if 
\begin{enumerate}
\item[(i)] $S_t: \Omega \to X_\partial$ is $\mathcal M/\mathscr B(X_\partial)$-measurable;
\item[(ii)] $\zeta: \Omega \to [0,\infty]$ is $\mathcal M$-measurable;
\item[(iii)] for any $\omega \in \Omega$, $S_t(\omega) \in X$ whenever $t<\zeta(\omega)$ and $S_t=\partial$ for all $t \ge \zeta(\omega)$.
\end{enumerate}
We say that $(\E, \F)$ is {\it local} if $\E(f,g)=0$ whenever $f, g \in \mathcal F$ with ${\rm supp}[f] \cap {\rm supp}[g]=\emptyset$.
We say that $ (\E, \F)$ is {\it strongly local} if for any $f,g \in \mathcal F$, the following holds: if $g$ is constant on a neighbourhood  of ${\rm supp}[f]$, then  $\E(f,g)=0.$
If a quasi-regular Dirichlet form $(\E,\F)$ is local, then the corresponding processes have continuous paths on $[0,\zeta)$ $\mathbb P^x$-almost surely for $\E$-q.e.\ $x \in X$ (see \cite[Theorem 1.11]{MR92}).

\subsection{RCD Spaces} \label{subsec: RCDLV}
In this subsection, we recall RCD$(K,\infty)$/RCD$^*(K,N)$ spaces. 
Recall that ${\sf Ch}$ denotes the Cheeger energy and the property of infinitesimal Hilbertianity (IH) was defined in \eqref{defn: Cheeger}.
Under (IH), ${\sf Ch}$ becomes a strongly local symmetric Dirichlet form (\cite{AGS14, AGS14b}). By the third paragraph in Subsection \ref{subsec: NSD}, there exists the corresponding semigroup $\{H_t\}_{ t \ge 0}$ (called {\it heat semigroup}) and the infinitesimal generator $\Delta$. 
Let us consider the following condition:
\begin{align} \label{text: REP}
&\text{Every function $f \in W^{1,2}(X,d,m)$ with $|\nabla f| \le 1$ $m$-a.e., admits a continuous}  \notag
\\
&\text{$1$-Lipschitz representative $\tilde{f}$.} 
\end{align}
We recall gradient estimates of the heat semigroups: for every $f \in W^{1,2}(X,d,m)$ with $|\nabla f| \le 1$ $m$-a.e., and every $t>0$, we have 
\begin{align} \label{BEI}
{H}_t f \in {\rm Lip}_b(X), \quad |\nabla {H}_tf|^2 \le e^{-2Kt}{H}_t(|\nabla f|^2), \quad m\text{-a.e.\ in $X$}.
\end{align}
The gradient estimate with dimensional upper bounds is as follows:
\begin{align}\label{BEK}
{H}_t f \in {\rm Lip}_b(X), \quad |\nabla {H}_tf|^2 +\frac{4Kt^2}{N(e^{2Kt}-1)}|\Delta H_t f|^2\le e^{-2Kt}{H}_t(|\nabla f|^2), \quad m\text{-a.e.\ in $X$}.
\end{align}
According to a sequence of results \cite{AGS15, AGS14, EKS15}, RCD$(K,\infty)$/RCD$^*(K,N)$ conditions can be identified with  \eqref{BEI}/\eqref{BEK} if we assume \eqref{asmp: basic}, \eqref{asmp: basic1}, (IH) and \eqref{text: REP}. 

A metric measure space $(X,d,m)$ is called {\it an RCD$(K,\infty)$ (resp.\ RCD$^*(K,N)$) space} if \eqref{BEI} (resp.\ \eqref{BEK}) holds under the assumptions \eqref{asmp: basic}, \eqref{asmp: basic1}, (IH) and \eqref{text: REP}.
\begin{rem} \normalfont
Note that RCD$(K,\infty)$/RCD$^*(K,N)$ conditions have been originally defined in terms of $K$-convexity of the relative entropy in the $L^2$-Wasserstein space with condition (IH) (\cite{AGMR15, AGS14, EKS15})). 
\end{rem}
The class of RCD spaces contains various  (finite- and infinite-dimensional) singular spaces such as Ricci limit spaces (Sturm \cite{Sturm06, Sturm06-2}, Lott--Villani \cite{LV09}), Alexandrov spaces (Petrunin, Zhang--Zhu \cite{Pet11, ZZ10}), warped products and cones (Ketterer \cite{Ket14, Ket14a}), quotient spaces (Galaz-Garc\'ia--Kell--Mondino--Sosa \cite{GKMS17}) and infinite-dimensional spaces such as Hilbert spaces with log-concave measures (Ambrosio--Savar\'e--Zambotti \cite{ASZ09}) (related to various stochastic partial differential equations). See these references for concrete examples. 

An important property of RCD$(K,\infty)$/RCD$^*(K,N)$ spaces is their stability under the pmG convergence. 
\begin{thm} [\cite{ AGMR15, AGS14b, EKS15, GMS13, Sturm06, Sturm06-2}] {\bf (Stability of the RCD$(K,\infty)$/RCD$^*(K,N)$)}
\\
Let $\mathcal X_n$ be an RCD$(K,\infty)$/RCD$^*(K,N)$ space for $n \in \N$.
If $\mathcal X_n$ converges to $\mathcal X_\infty$ in the pmG sense, then the limit space $\mathcal X_\infty$ also satisfies the RCD$(K,\infty)$/RCD$^*(K,N)$ condition. 
\end{thm}

\section{Construction of Non-Symmetric Dirichlet Forms} \label{sec: NSDF}
In this section, we construct a non-symmetric Dirichlet form consisting of a uniformly elliptic operator,\deleted{a} derivations and a killing term.\
Let us consider the following bilinear form $\mathcal E: {\rm Lip}_{bs}(X) \times {\rm Lip}_{bs}(X) \to \R$:
$$\mathcal E (f,g):=\frac{1}{2}\int_X\langle A\nabla f, \nabla g\rangle dm+\int_{X}\b_1(f)gdm+\int_{X}f\b_2(g)dm+\int_Xfgcdm.$$
We write $\int_{X}\langle A \nabla f, \nabla g \rangle dm=2{\sf Ch}_A(f,g)$ in short.
Recall that $L^2(TX)$ denotes the tangent module as in Definition \ref{defn: tan}.\ 
Let $A: L^{2}(TX) \to L^2(TX)$ denote a (not necessarily $L^2(TM)$-symmetric) \replaced{module morphism}{linear operator} satisfying that 
there exists $H \in L^1_{loc}(X,m)$ so that $|A {Y}| \le H |{Y}|$ for any {$Y \in L^2(TX)$}.\ We write $|A|$ for the minimal element among such $H$.\
Let $\tilde{A}:=1/2(A+A^*)$ and $\check{A}:=1/2(A-A^*)$ where $A^*$ is the $L^2(TX)$-adjoint operator of $A$. 
Suppose the following conditions:
\begin{asmp}\label{asmp: NS1} \normalfont
Let $(X,d,m)$ be a metric measure space with \eqref{asmp: basic}, \eqref{asmp: basic1} and (IH). 
Suppose that
\begin{itemize}
\item[(i)] $|\check{A}| \in L^\infty(X,m)$ and there exists a constant $\lambda>0$ so that for any $f \in {\rm Lip}_{bs}(X)$, 
\begin{align} \label{ineq: UELB}
\langle A \nabla f, \nabla f \rangle \ge \lambda \langle \nabla f, \nabla f \rangle, \ m\text{-a.e.};
\end{align}
\item[(ii)] ${|\b_1+\b_2| , c \in L^\infty_{loc}(X,m)}, |\b_1-\b_2| \in L^\infty(X,m)$;
\item[(iii)] for any non-negative $f \in {\rm Lip}_{bs}(X), i=1,2,$
\begin{align} \label{ineq: DSM}
\int_X(\b_i(f)+cf)dm \ge 0.
\end{align}
\end{itemize}
\end{asmp}

Then we can show the following proposition.
\begin{prop} \label{thm: NS1}
Suppose Assumption \ref{asmp: NS1}.
\begin{description}
\item[(a)] If $(X,d,m)$ is an RCD$^*(K,N)$ space, then the form \eqref{form: non-sym} is closable and the closed form $(\mathcal E, \F)$ is a regular local Dirichlet form.
\item[(b)]  Let $(X,d,m)$ be an RCD$(K,\infty)$ space with $m(X)<\infty$. If $|A|, |\b_1|, |\b_2|, c \in L^\infty(X,m)$,  
then the form \eqref{form: non-sym} is closable and the closed form $(\mathcal E, \F)$ is a quasi-regular local Dirichlet form.
\end{description}
\end{prop}

{\it Proof of Proposition \ref{thm: NS1}}. 
{{Non-negativity $(\E.1)$}:} 
{By $|A| \in L^1_{loc}(X,m)$, the integrand $\langle A \nabla f, \nabla g \rangle$ for $f, g \in {\rm Lip}_{bs}(X)$ is $m$-integrable.}
For $f \in {\rm Lip}_{bs}(X)$, by Leibniz formula of $\b_i$, \eqref{ineq: UELB} and \eqref{ineq: DSM}, we have
\begin{align} \label{eq:LBB}
\E(f) = {\sf Ch}_A(f)+\int_{X}\Bigl(\b_1(f^2)+\b_2(f^2)+2f^2c\Bigr)dm \ge  \lambda{\sf Ch}(f)\ge 0.
\end{align}
\deleted{Thus we see that ${\E_{{\alpha_0}}(f)=}\E(f){+\alpha_0\|f\|^2_2} \ge \lambda/2 {\sf Ch}(f) \ge 0$.}

{Closability:} Let $\{f_n\}_{n \in \N} \subset \L_{bs}(X)$ be an $\mathcal E$-Cauchy sequence so that $\|f_n\|_2 \to 0$.
Since $\C$ is closable by \cite{AGS14, AGS14b}, we have $\C(f_n) \to 0$, which implies that there exists a subsequence $f_{n'}$ so that $|\nabla f_{n'}|$ converges to zero $m$-a.e.. By $|A\nabla f_n| \le |A||\nabla f_n|$, we have that $|A\nabla f_{n'}|\to 0$ $m$-a.e. Noting that ${\sf Ch}_A(f) \le 2\E(f)$ by \eqref{eq:LBB}, and by using Fatou's lemma,  we have 
\begin{align*}
&2{\sf Ch}_A(f_n)=\int_{X}\lim_{n' \to \infty}\Bigl\langle A(\nabla f_n-\nabla f_{n'}), \nabla f_n-\nabla f_{n'} \Bigr\rangle dm 
\\
&\le \liminf_{n' \to \infty}\int_{X}\Bigl\langle A\nabla (f_n-f_{n'}), \nabla (f_n-f_{n'}) \Bigr\rangle dm  \le \liminf_{n' \to \infty}\E(f_n-f_{n'}).
\end{align*}
Since $\{f_n\}_{n \in \N}$ is an $\E$-Cauchy sequence, the R.H.S.\ above can be arbitrarily close to zero as $n$ is sufficiently large. Thus, \replaced{${\sf Ch}_A(f_n) \to 0$ as $n \to \infty$.}{${\sf Ch}_A$ is closable.} 
We next show the closability of the remaining part of $\E$. 
By replacing $f_n$ to $f_nh$ for any $h \in {\rm Lip}_{bs}(X)$ with $h \ge 0$, we may assume that ${\rm supp}(f_n) \subset K \subset X$ for some bounded open set $K$.
Noting that  $\b_i(f_{n}) \le |\b_i||\nabla f_{n}|$ for $i=1,2$ and ${\sf Ch}(f_n) \to 0$ by the closability of ${\sf Ch}$, we have
\begin{align*}
&\Bigl|\int_{X}\b_1(f_n)f_ndm+\int_{X}f_n\b_2(f_n)dm+\int_Xf_n^2cdm \Bigr| 
\\
&\le\bigl\||\b_1+\b_2|f_n\1_{K}\bigr\|^{1/2}_2{\sf Ch}(f_n)^{1/2}+\|c\1_Kf^2_n\|_1
\\
&\le \bigl\||\b_1+\b_2|\1_{K}\bigr\|^{1/2}_\infty\|f_n\|^{1/2}_2{\sf Ch}(f_n)^{1/2} +\|c\1_K\|_\infty \|f_n\|_2
\\
& \to 0 \quad \text{($n \to \infty$)}.
\end{align*}
Thus, we have $\E(f_n) \to 0$ and we have proved the closability of $\E$. 

{Weak sector condition $(\E.2)$}: \deleted{We first give a proof for the case of $\alpha_0=0$. In this case, }\replaced{I}{i}t suffices to show that there exists a constant ${C}>0$ so that $\check{\E}(f,g) \le {C}\E_{1}(f)\E_{{1}}(g)$, whereby $\check{\E}(f,g)=1/2(\E(f,g)-\E(g,f))$ (see \cite[Chapter I, \S 2]{MR92}). \deleted{Let $B=B_r(\x) \supset {\rm supp}[f]\cup {\rm supp}[g]$ be an open ball containing the support of $f,g \in \L_{bs}(X)$.} Then
\begin{align}
2\check{\E}(f,g)&=\int_{X}\Bigl(\frac{1}{2}\langle \check{A}\nabla f, \nabla g\rangle+\b_1(f)g-\b_2(f)g-f\b_1(g)+f\b_2(g) \Bigr)dm \notag
\\
&\le \frac{1}{2}\bigl\||\check{A}|\bigr\|_\infty\|\nabla f\|_2\|\nabla g\|_2+\|\b_1(f)-\b_2(f)\|_2\|g\|_2+\|\b_2(g)-\b_1(g)\|_2\|f\|_2 \notag
\\
&\le  \frac{1}{2}\bigl\||\check{A}|\bigr\|_\infty\|\nabla f\|_2\|\nabla g\|_2+\bigl\| |\b_1-\b_2||\nabla f|\bigr\|_2\|g\|_2+ \bigl\||\b_1-\b_2||\nabla g|\bigr\|_2\|f\|_2.
\label{formula: SEC}
\end{align}
Since $|\b_1-\b_2| \in L^\infty(X;m)$ and $2{\sf Ch}(f)=\|\nabla f\|_2^2 \le (2/\lambda)\E(f)$ by \eqref{eq:LBB}, we have 
\begin{align*}
&\frac{1}{2}\bigl\||\check{A}|\bigr\|_\infty\|\nabla f\|_2\|\nabla g\|_2+\bigl\| |\b_1-\b_2||\nabla f|\bigr\|_2\|g\|_2+ \bigl\| |\b_1-\b_2||\nabla g|\bigr\|_2\|f\|_2
\\
&\le  \frac{1}{\lambda}\bigl\||\check{A}|\bigr\|_\infty \E_{1}^{1/2}(f)\E_{1}^{1/2}(g)+  \sqrt{\frac{2}{\lambda}}\|(\b_1-\b_2)\|_\infty \E_{1}^{1/2}(f)\E_{1}^{1/2}(g)
\\
&\quad + \sqrt{\frac{2}{\lambda}}\|(\b_1-\b_2)\|_\infty\E_{1}^{1/2}(g)\E_{1}^{1/2}(f)
\\
& \le {C}\E_{1}^{1/2}(f)\E_{1}^{1/2}(g).
\end{align*}
\deleted{Taking $\alpha_0=0$,} \replaced{W}{w}e finish\added{ed} to prove that $\E$ satisfies the weak sector condition.

{\replaced{Markov Property $(\E.4)$}{Markov property $(\E.4)$ and Dual Sub-Markov Property $(\hat{\E}.4)$}}: 
Let $\phi_\e: \R \to [-\e, 1+\e]$ be an infinitely differentiable function so that $0 \le \phi_\e(t)-\phi_\e(s) \le t-s$ for all 
$t,s \in \R$ with $t \ge s$ and
\begin{align*}
\phi_\e(t)=
\begin{cases}
t \quad \hspace{7mm}\text{for $t \in [0,1]$},
\\
1+\e \quad \text{for $t \in [1+2\e, \infty)$},
\\
-\e\quad \hspace{4mm} \text{for $t \in (-\infty, -2\e]$}.
\end{cases}
\end{align*}
By \cite[Proposition 4.7, 4.10 in Chapter I]{MR92}, 
it suffices to show that for all $f \in {\rm Lip}_{bs}(X)$ and $\e>0$, it holds that $\phi_\e(f) \in \mathcal F$ and 
\begin{align} \label{ineq: SUB-DUAL}
\liminf_{\e \to 0}\E(\phi_\e(f), f-\phi_\e(f)) \ge 0, \quad \liminf_{\e \to 0}\E(f-\phi_\e(f),\phi_\e(f)) \ge 0.
\end{align}
It is clear that $\phi_\e(f) \in {\rm Lip}_{bs}(X) \subset \mathcal F$.
We only show the L.H.S.\ side of \eqref{ineq: SUB-DUAL} (the proof of the R.H.S.\ is similar). The diffusion part ${\sf Ch}_A$ clearly satisfies \eqref{ineq: SUB-DUAL} by the chain rule, $0 \le \phi'_\e \le 1$, \eqref{prop: mod}  and \eqref{eq:LBB}. In fact, 
$$\int_{X}\Bigl\langle A \nabla \phi_\e(f), \nabla \bigl(f- \phi_\e(f)\bigr)\Bigr\rangle dm=\int\phi'_\e(f)(1-\phi'_\e(f))\langle A \nabla f, \nabla f \rangle dm \ge 0.$$
\added{For the remaining part,} we have
\begin{align}
&\int_X\b_1(f-\phi_\e(f))(\phi_\e(f))dm+\int_X(f-\phi_\e(f))\b_2(\phi_\e(f))dm + \int_X(f-\phi_\e(f))(\phi_\e(f))cdm \notag
\\
& = \int_X\Bigl(\b_1\Bigl(\phi_\e(f)(f-\phi_\e(f))\Bigr)-\phi_\e(f)(f-\phi_\e(f))c\Bigr) dm+\int_X\phi'_\e(f-\phi_\e(f))(\b_1-\b_2)(f)dm. \notag
\end{align}
The first term in the second line above is non-negative since \eqref{ineq: DSM} and $(f-\phi_\e(f))\phi_\e(f) \ge 0$. The second term converges to zero since 
$$\phi'_\e(f)(f-\phi_\e(f)) \le (\1_{[-2\e, 1+2\e]}(f))(f-\phi_\e(f)) \to 0 \quad \e \downarrow 0.$$
Thus, we finished to prove the Markovian property $(\E.4)$.

\deleted{In the case of $\alpha_0=0$, the last integral is non-negative. 
The dual form $\hat{\E}(f,g)=\E(g,f)$ also satisfies the Markov property by the same argument.
In the case of $\alpha_0>0$ with $\b_2=0$, the last integral is clearly non-negative. In the case of $\alpha_0>0$ with $\b_1=0$, the dual form satisfies the Markov property by the same argument. thus we have finished the proof of the Markov property. The last integral is non-negative. Thus we have proved the sub-Markov property of $\E$. The sub-Markov property of the dual form $\hat{\E}$ can be proved by the same manner.}
In the case of RCD$^*(K,N)$ spaces, the underlying space $X$ is locally compact by the local volume doubling property according to the Bishop--Gromov inequality \cite[Proposition 3.6]{EKS15} (see also \cite[{Corollary} 2.4]{Sturm06-2}). Therefore, $(\E, \F)$ is regular since ${\rm Lip}_{bs}(X)$ is dense in both in $C_0(X)$ with $\|\cdot\|_\infty$ and $\mathcal F$ with $\tilde{\E}_1^{1/2}$.
In the case of RCD$(K,\infty)$ spaces, \deleted{the state space}$X$ is generally not locally compact and therefore we need to show \added{that there exists an $\E$-nest of compact sets $\{E_k\}_{k \in \N}$, which is called} tightness of capacity.
\\
{Tightness of Capacity:} 
We show that there exists an $\E$-nest of compact sets $\{E_k\}_{k \in \N}$, i.e., an increasing sequence of compact sets $E_k \subset E_{k+1}$ so that ${\rm Cap}_{\E}(X \setminus E_k) \to 0$. Let $\{x_k\}_{k \in \N} \subset X$ be a countable dense subset.
Define 
$$w_n(x):=\min\{1, \min_{1 \le k \le n}d(x_k, x)\}.$$
We see that $0 \le w_n \le 1$ and $w_n \downarrow 0$ as $n \to \infty$. Thus, $w_n \to 0$ in $L^2(m)$. By the definition of ${\sf Ch}$, it is easy to see that $2\C(w_n) \le  {\rm Lip}(w_n) \le 1.$ 
Noting $|\nabla w_n| \le 1$, we have that 
\begin{align} \label{Quasi-Reg}
{\E}(w_n)
&=\C_{{A}}(w_n)+\int_X(\b_1(w_n)w_n)dm+\int_X(w_n\b_2(w_n))dm+\int_Xw_n^2cdm \notag
\\
&\le \frac{\bigl\| |A| \bigr\|_\infty}{2}+ \||\b_1+\b_2|\|_\infty\|\nabla w_n\|_2\|w_n\|_2+\|c\|_\infty\|w_n\|^2_2 \notag
\\
&\le \frac{\bigl\| |A| \bigr\|_\infty}{2} + \||\b_1+\b_2|\|_\infty+\|c\|_\infty.
\end{align}
%
Therefore, $w_n$ is a uniformly bounded sequence in \replaced{$\mathcal F$}{$W^{1,2}(X,d,m)$} with respect to $\tilde{\mathcal E}_{{1}}$ where $\tilde{\E}_{{1}}(f,g):=1/2(\E_{{1}}(f,g)+\E_{{1}}(f,g)).$  Thus, $w_n \to 0$ weakly in \replaced{$\mathcal F$}{$W^{1,2}(X,d,m)$} with respect to $\tilde{\mathcal E}_{{1}}$ by \cite[Lemma 2.12 in Chapter I]{MR92}. By the Banach-Saks theorem, there exists a subsequence $\{n(i)\}_{i}$ so that the Cesaro means 
$$v_l=\frac{1}{l}\sum_{i=1}^lw_{n(i)}$$
converges\deleted{strongly} in \replaced{$\mathcal F$}{$W^{1,2}(X,d,m)$ with respect to $\tilde{\mathcal E}_{{1}}$}, thus in $(L^2(X,m), \|\cdot\|_2)$.  
By \cite[Proposition III. 3.5]{MR92}, we have that there exists a subsequence $\{v_{i(j)}\}_j$ so that $v_{i(j)} \to 0$ quasi-uniformly, i.e., for any $k$ there exists a closed set $G_k \subset X$ so that ${\rm Cap}_\E(X\setminus G_k) \le 1/k$ and $v_{i(j)} \to 0$ uniformly on $G_k$. Since $w_{n(i(j))} \le v_{i(j)}$, setting $F_k=\cap_{i \le k}G_i$, we have that $w_{n(i(j))} \to 0$ uniformly on $F_k$ for all $k$ and ${\rm Cap}_\E(X\setminus F_k) \le 1/k$.
Let $\e>0$ and $n$ be an integer so that $w_n <\e$, and the definition of $w_n$ implies
$$F_k \subset \bigcup_{k=1}^nB(x_k,\e).$$
Since $\e$ is arbitrarily small, we have that $F_k$ is totally bounded and thus compact. 
{The other conditions for quasi-regularity in \cite[(ii) and (iii) in Definition 3.1 in Chapter IV]{MR92} are easy to check since ${\rm Lip}_{bs}(X)$ separates points and is dense in $\F$ with respect to $\tilde{\E}^{{1/2}}_1$.}

The local property is obvious according to the locality of ${\sf Ch}$ and derivations $\b_i$ for $i=1,2$. Thus, we have finished the proof of Proposition \ref{thm: NS1}\deleted{under Assumption \ref{asmp: NS1}}. 
\qed

\begin{rem} \normalfont \label{rem: DA}
Construction of non-symmetric diffusion on metric measure spaces (including RCD spaces) has been already considered in Fitzsimmons \cite{F07} and Trevisan \cite{T14} with different approaches and different scopes. 
\begin{itemize}
\item[(i)] In the former paper, the Girsanov transform was used to produce drift perturbations from symmetric diffusions. 
\item[(ii)] In the latter paper, the martingale problem was developed to construct diffusion processes in this generality. 
\end{itemize}
An advantage of the Dirichlet form approach adopted in this paper is to make the issue of convergence \replaced{simpler}{easier compared to the other two approaches}. \added{This is mainly because the domains of Dirichlet forms corresponding to \eqref{form: non-sym} has the common core ${\rm Lip}_{bs}(X)$, which is useful especially to show tightness of diffusion processes in Section \ref{sec: CFDDNS}}. 
\deleted{The first reason is that the dual sub-Markov property (does not necessarily hold in the other two situations) enables us to utilize the Lyons-Zheng decomposition for tightness. The second reason is that we do not need to deal with varying domains of generators since we can take a common core ${\rm Lip}_{bs}(X)$ for $(\E_n, \F_n)$ after embedding $X_n$ into $X$, which helps to show the tightness of processes.} 
\end{rem}

\section{Convergence of Non-symmetric Forms} \label{sec: CNSF}
In this section, we show the convergence of non-symmetric forms. We first modify the definition in \cite{H98} for varying metric measure spaces and prove this modified convergence under Assumption \ref{asmp: NS3}. We recall the definition of the $L^p$-convergence on varying metric measure spaces in the sense of pmG following \cite{GMS13}.

\begin{defn}$(${\rm See \cite[Definition 6.1]{GMS13}}$)$ \label{defn: Weak} \normalfont
Let $(X_n, d_n, m_n, {\x_n})$ be a sequence of p.m.m for $n \in \EN$. spaces. Assume that $(X_n,d_n,m_n, \x_n)$ converges to $(X_\infty, d_\infty, m_\infty, \x_\infty)$ in the pmG sense.  Let $(X,d)$ be a complete separable metric space and
$\iota_n: \mathrm{supp}[m_n] \to X$ be isometries as in Definition \ref{prop: Dconv}. We identify $(X_n, d_n, m_n)$ with $(\iota_n(X_n), d, {\iota_n}_\#m_n)$ and omit $\iota_n$. 
\begin{description}
	\item[(i)] We say that {\it $f_n \in L^2(X, m_n)$ converges weakly to $f_\infty \in L^2(X,m_\infty)$} if the following hold:
$$\sup_{n \in \N}\int|f_n	|^2\ dm_n < \infty \quad \text{and} \quad \int \phi f_n \ dm_n \to \int \phi f_\infty \ dm_\infty \quad \forall \phi \in C_{bs}(X).$$
 \item[(ii)] We say that {\it $f_n \in L^2(X,m_n)$ converges strongly to $f_\infty \in L^2(X, m_\infty)$} if $f_n$ converges weakly to $f_\infty$ and the following holds:
$$\limsup_{n \to \infty}\int |f_n|^2 \ dm_n \le \int|f_\infty|^2\ dm_\infty.$$
\end{description}
\end{defn}

Now we introduce a modified definition of the convergence of non-symmetric forms in \cite{H98} for varying metric measure spaces in the sense of the pmG.
The modified point is just to replace the usual $L^2$-convergence with the $L^2$-convergence in Definition \ref{defn: Weak}.
Let $\mathcal X_n=(X_n, d_n, m_n, \x_n)$ be a sequence of p.m.m. spaces converging to a limit $\mathcal X_\infty=(X_\infty, d_\infty, m_\infty, \x_\infty)$ in the pmG sense. Let $(X,d)$ be a complete separable metric space and
$\iota_n: \mathrm{supp}[m_n] \to X$ be isometries as in Definition \ref{prop: Dconv}. Let $(\E_n,\F_n)$ be a sequence of coercive forms on $L^2(X;m_n)$. 

Let $\Phi_n(f):=\sup\{\E^n_{1}(g, f): \tilde{\E}^n_{1}(g)^{1/2}=1\}$ for $f \in \mathcal F_n$, where\added{by} $\tilde{\E}^n_{1}$ denotes the symmetric part of $\E^n_{1}$.
\begin{defn} \label{thm: Mosco of Ch} \normalfont (See also \cite[Definition 7.11]{T10})
 We say that {\it $(\E_n, \F_n)$ converges to $(\E_\infty, \F_\infty)$} if the following two conditions hold:
\begin{description}
	\item[(F1)] If a sequence $f_n \in L^2(X,m_n)$ converges weakly to $f_\infty \in L^2(X,m_\infty)$ with $\liminf_{n \to \infty}\Phi_n(f_n)<\infty$, then 
	it holds that
	$$f_\infty \in \mathcal F_\infty.$$
	\item[(F2)] For any sequence $f_n \in \F_n$ converging weakly in $L^2$ to $f_\infty \in \F_\infty$, and any $w_\infty \in \mathcal F_\infty$, there exists a sequence $w_n \in \mathcal F_n$ converging strongly in $L^2$ to $w_\infty \in \mathcal F_\infty$ so that 
	$$\lim_{n \to \infty}\mathcal E_n(f_n,w_n)=\mathcal E_\infty(f_\infty, w_\infty).$$
\end{description}
\end{defn}
\begin{rem} \normalfont
In T\"olle \cite{T06, T10}, he introduced a notion of a convergence of non-symmetric forms whose basic Banach spaces vary. The difference of his approach and this paper is the notion of the $L^2$-convergence on varying metric measure spaces whereby in this paper we follow \cite{GMS13}. If \replaced{the}{our} Hilbert spaces $\{L^2(X;m_n)\}_{n \in \EN}$ have an asymptotic relation in \cite{T10}, these two notions of the $L^2$-convergence are equivalent, so the following Theorem \ref{thm: EQUIVHINO} corresponds to \cite[Theorem 7.15]{T10}.
\end{rem}
Verifying (F2) is not always easy and we introduce another condition\deleted{ for this purpose}:
\begin{description}
	\item[(F2$'$)] For any sequence $\{n_{{k}}\} \uparrow \infty$ and any sequence $f_k \in L^2(X; m_{n_k})$ weakly convergent in $L^2$ to $f_\infty \in \mathcal F_\infty$ with $\sup_{k \in \N}\Phi_{n_k}(f_k)<\infty$, there exists a dense subset $\mathcal C \subset \mathcal F_\infty$ for the topology with respect to $\tilde{\E}^\infty_{1}$ so that every $w \in \mathcal C$ has a sequence $\{w_k\}$ \added{with $w_k \in \mathcal F_{n_k}$} converging to $w$ strongly \added{in $L^2$} with 
	$$\liminf_{k \to \infty}\E_{n_k}(w_k,f_k) \le \E_\infty(w,f_\infty).$$
\end{description}
We also define (F$1_*$) by replacing $\Phi_n(f_n)$ with $\tilde{\E}^n_1(f_n)^{1/2}$ in (F1), and (F$2_*'$) by replacing $\Phi_{n_k}(f_k)$ with $\tilde{\E}^{n_k}_1(f_k)^{1/2}$ in (F2$'$).

\replaced{We now study}{What we would like to study in this section is to investigate} the relation\deleted{ship} between the convergence of forms and $L^2$-convergences of the corresponding semigroups and resolvents. Let $\{T^n_t\}_{t\ge 0}$ and $\{G_\alpha^n\}_{\alpha \ge 0}$ be the $L^2$-contraction semigroup and resolvent associated with $\mathcal E_n$. 
\begin{description}
	\item[(R)] For any sequence ${f_n}$ converging to ${f_\infty}$ strongly in $L^2$, the resolvent $G^n_{{\alpha}} f_n$ converges to $G_{{\alpha}}^\infty f_\infty$ strongly in $L^2$ for any ${\alpha>0}$.
	\item[(S)] For any sequence ${f_n}$ converging to ${f_\infty}$ strongly in $L^2$, $T^n_t f_n$ converges to $T_t^\infty f_\infty$ strongly in $L^2$. The convergence is uniform on any compact time interval $[0,T].$ 
\end{description}

\begin{thm}  \label{thm: EQUIVHINO}
The following statements hold:
\begin{enumerate}
\item[(i)] $(F1) (F2) \iff (F1)(F2') \iff (R) \iff (S)$;
\item[(ii)]  $(F1_*)(F2_*') \implies (R)$.
\end{enumerate}
\end{thm}
\proof
The proof is just a modification of \cite[Theorem 3.2, Corollary 3.3]{H98} for varying metric measure spaces in the sense of the pmG so that the notion of $L^2$-convergence is replaced by Definition \ref{defn: Weak}. So we omit the proof. 
\qed

Hereafter in this section, we focus on the convergence of Dirichlet forms corresponding to $\eqref{form: non-sym}$.
To characterize the convergence of these forms in terms of convergences coefficient\added{s}, we introduce a convergence of $A_n$ and recall a convergence of derivation operators $\b_i^n$ {(\cite{AST16})}.
Let $\mathscr A\subset {\rm Lip}_b(X_\infty)$ denote the smallest algebra containing the following functions:
$$\min\{d(\cdot, x), k\},\quad k \in \Q, \ x \in D \subset X_\infty, \ \text{dense subset}.$$
The algebra $\mathscr A$ becomes a vector space over $\Q$. Let $\mathscr A_{bs}$ be a subalgebra consisting of bounded support functions. 
Let $\{H^\infty_t\}_{t \ge 0}$ be the heat semigroup associated with Cheeger energy ${\sf Ch}_\infty$ (note that $\{H^\infty_t\}_{t \ge 0}$ is not the semigroup associated with the non-symmetric form $\mathcal E$).
Let $H_{\Q_+}\mathscr A_{bs}:=\{H^\infty_sf: f \in \mathscr A_{bs}, s \in \Q_+\} \subset {\rm Lip}_{b}(X)$, where\added{by} ${\rm Lip}_{b}(X)$ denotes the set of bounded Lipschitz functions on $X$.
\replaced{Recall that}{Let} ${\rm Der}^p_{loc}(X,d,m)$ and ${\rm Der}^p(X,d,m)$ be the set of derivation operators $\b$ in $(X,d,m)$ with $|\b| \in L^p_{loc}(X,d,m)$ and $|\b| \in L^p(X,d,m)$, respectively. 
\begin{defn} \normalfont $(${\rm \cite[Definition 4.3, 5.3]{AST16}}$)$ Let $(X_n,d_n,m_n,\x_n)$ converge to $(X_\infty, d_\infty, m_\infty, \x_\infty)$ in the pmG sense. Let $(X,d)$ be a complete separable metric space and
$\iota_n: \mathrm{supp}[m_n] \to X$ be isometries as in Definition \ref{prop: Dconv}.
\begin{itemize}
\item[(i)] (Weak Convergence) We say that $\b_n \in {\rm Der}^1_{loc}(X,d,m_n)$ converges weakly to $\b_\infty \in  {\rm Der}^1_{loc}(X,d,m_\infty)$  in duality with $H_{\Q_+}\mathscr A_{bs}$ if, for all $f \in H_{\Q_+}\mathscr A_{bs}$, 
$$\int_X\b_n(f)h dm_n \to \int_X\b_\infty(f)h dm_\infty \quad \forall h \in C_{bs}(X).$$
\item[(ii)] (Strong Convergence) We say that $\b_n \in {\rm Der}^1_{loc}(X,d,m_n)$ converges strongly to $\b_\infty \in  {\rm Der}^1_{loc}(X,d,m_\infty)$ if, for all $f \in H_{\Q_+}\mathscr A_{bs}$, the function $\b_n(f)$ converges in measure to $\b_\infty(f)$, i.e., 
$$\int_X\Phi(\b_n(f))h dm_n \to \int_X\Phi(\b_\infty(f))h dm_\infty \quad \forall h \in C_{bs}(X), \quad \forall \Phi \in C_b(\R).$$
\item[(iii)] ($L^p$-strong Convergence) Let $p \in [1,\infty)$. We say that $\b_n \in {\rm Der}^p_{loc}(X,d,m_n)$ converges $L^p_{loc}$-strongly to $\b_\infty \in  {\rm Der}^p_{loc}(X,d,m_\infty)$  if $\b_n$ converges strongly to $\b_\infty$ and for all $f\in H_{\Q_+}\mathscr A_{bs}$ and $R>0$,  
$$\limsup_{n \to \infty} \int_{B_R(\x_n)}|\b_n(f)|^pdm_n \le \int_{B_R(\x_\infty)}|\b_\infty(f)|^pdm_\infty.$$
Analogously we say that $\b_n \in {\rm Der}^p_{loc}(X,d,m_n)$ converges $L^p$-strongly to $\b_\infty \in  {\rm Der}^p_{loc}(X,d,m_\infty)$  if, $\b_n$ converges strongly to $\b_\infty$ and, for all $f\in H_{\Q_+}\mathscr A_{bs}$,
$$\limsup_{n \to \infty} \int_{X}|\b_n(f)|^pdm_n \le \int_{X}|\b_\infty(f)|^pdm_\infty.$$
\end{itemize}
\end{defn}

Now we recall the $W^{1,2}$-convergence of functions on varying metric measure spaces in the sense of pmG. 
\begin{defn}$(${\rm\cite[Definition 5.2]{AH16}}$)$ \label{defn: WConv} \normalfont
Let $(X_n, d_n, m_n, \x_n)$ be a sequence of p.m.m. spaces. Assume that $(X_n,d_n,m_n, \x_n)$ converges to $(X_\infty, d_\infty, m_\infty, \x_\infty)$ in the pmG sense.  Let $(X,d)$ be a complete separable metric space and
$\iota_n: \mathrm{supp}[m_n] \to X$ be isometries as in Definition \ref{prop: Dconv}. We identify $(X_n, d_n, m_n)$ with $(\iota_n(X_n), d, {\iota_n}_\#m_n)$ and omit $\iota_n$. 
\begin{description}
	\item[(i)] We say that {\it $f_n \in W^{1,2}(X, m_n)$ converges weakly to $f_\infty \in W^{1,2}(X,m_\infty)$ in $W^{1,2}$} if $f_n \to f_\infty$ weakly in $L^2$ in the sense of Definition \ref{defn: Weak} and 
$\sup_{n \in \N}{\sf Ch}_n(f_n)< \infty$;
 \item[(ii)] We say that {\it $f_n \in W^{1,2}(X,m_n)$ converges strongly to $f_\infty \in W^{1,2}(X, m_\infty)$ in $W^{1,2}$} if $f_n$ converges strongly to $f_\infty$ in $L^2$ in the sense of Definition \ref{defn: Weak} and $\lim_{n \to \infty}{\sf Ch}_n(f_n)={\sf Ch}_\infty(f_\infty)$.
\end{description}
\end{defn}

Now we introduce a convergence of $A_n$.
\begin{defn} \normalfont \label{defn: CONVAN}
We say that $A_n$ {\it converges to} $A_\infty$ if for any $u_n \to u_\infty$ weakly in $W^{1,2}$ and $v_n \to v_\infty$ strongly in $W^{1,2}$, 
\begin{align*}
&\int_{X}\langle \nabla u_n, A_{n}\nabla v_n \rangle dm_{n} \to \int_{X}\langle \nabla u_\infty, A_\infty\nabla v_\infty \rangle dm_\infty,
\\
&\int_{X}\langle A_n\nabla u_n, \nabla v_n \rangle dm_{n} \to \int_{X}\langle A_\infty\nabla u_\infty, \nabla v_\infty \rangle dm_\infty.
\end{align*}
\end{defn}
Now we show the main theorem in this section.
\begin{thm}  \label{thm: HINOC}
Under Assumption \ref{asmp: NS3}, $(\mathcal E_n, \F_n)$ $($resp.\ $(\hat{\E}_n, \F_n)$$)$ converges to $(\mathcal E_\infty, \F_\infty)$ $($resp.\ $(\hat{\E}_\infty, \F_\infty)$$)$.
\end{thm}
\proof
By Theorem \ref{thm: EQUIVHINO}, it suffices to show $(F1_*)$ and $(F2_*')$.

$(F1_*)$: Let $u_n \to u_\infty$ weakly in $L^2$ and we may assume $\liminf_{n\to \infty}\mathcal E_n(u_n)<\infty$. Since $({\sf Ch}_n, \F_n)$ converges to $({\sf Ch}_\infty, \F_\infty)$ in the Mosco sense \cite[Theorem 6.8]{GMS13}, we have, by \eqref{eq:LBB}, 
\begin{align} \label{ineq: lambdaud}
{\sf Ch}_\infty(u_\infty) \le \liminf_{n \to \infty}{\sf Ch}_n(u_n)\le \frac{1}{\lambda}\liminf_{n\to \infty}\mathcal E^n(u_n)<\infty.
\end{align}
 This implies $u_\infty \in W^{1,2}(m_\infty)$ by the definition of $W^{1,2}(m_\infty)$.

$(F2_*')$: 
Let $n_k \uparrow \infty$ and $u_k \to u_\infty$ weakly in $L^2$ with $\sup_{k \in \N}\mathcal E_{n_k}(u_{k})<\infty$ and $u_\infty \in W^{1,2}(m_\infty)$. Then we have that $u_k \to u_\infty$ weakly in $W^{1,2}$ by definition. Let us take $\mathcal C=H_{\Q_+}\mathscr A_{bs}.$ Take $w \in \mathcal C$. 
By \eqref{eq: DIV}, we have 
\begin{align*}
&|\mathcal E_{n_k}(u_k, w)-\mathcal E_{\infty}(u,w)|
\\
&=\Bigl|{\sf Ch}^{n_k}_A(u_k,w)-{\sf Ch}_A^\infty(u_\infty,w) \Bigr|+\Bigl|\int_{X_{n_k}}\b^{n_k}_1(u_k)wdm_{n_k}-\int_{X_\infty}\b^{\infty}_1(u_\infty)wdm_\infty\Bigr|
\\
& \quad + \Bigl|\int_{X_{n_k}}u_k\b^{n_k}_2(w)dm_{n_k}-\int_{X_\infty}u_\infty\b^{\infty}_2(w)dm_\infty\Bigr|+\Bigl|\int_{X_{n_k}}u_kwc_{n_k}dm_{n_k}-\int_{X_{\infty}}u_\infty wc_{\infty}dm_{\infty}\Bigr|
\\
&=\Bigl|{\sf Ch}^{n_k}_A(u_k,w)-{\sf Ch}^\infty_A(u_\infty,w) \Bigr|+\Bigl|\int_{X_{n_k}}u_k\b^{n_k}_1(w)dm_{n_k}-\int_{X_\infty}u_\infty\b^{\infty}_1(w)dm_\infty\Bigr|\\
&\quad +
\Bigl|\int_{X_{n_k}}u_kw{\rm div} \b^{n_k}_1dm_{n_k}-\int_{X_\infty}u_\infty w {\rm div} \b^{\infty}_1 dm_\infty\Bigr|+ \Bigl|\int_{X_{n_k}}u_k\b^{n_k}_2(w)dm_{n_k}-\int_{X_\infty}u_\infty\b^{\infty}_2(w)dm_\infty\Bigr|
\\
& \quad +\Bigl|\int_{X_{n_k}}u_kwc_{n_k}dm_{n_k}-\int_{X_{\infty}}u_\infty wc_{\infty}dm_{\infty}\Bigr|
\\
&:={\rm (I)}_k+{\rm (II)}_k+{\rm (III)}_k+{\rm (IV)}_k+{\rm (V)}_k.
\end{align*} 
We first show ${\rm (I)}_k \to 0$ as $k \to \infty.$
By the convergence of $A_n$ to $A_\infty$, we have 
\begin{align*}
\Bigl|{\sf Ch}^{n_k}_A(u_k, w)-{\sf Ch}_A^\infty(u_\infty, w) \Bigr|=\frac{1}{2}\Bigl| \int_{X}\langle  A_{n_k}\nabla u_k,\nabla w \rangle dm_{n_k}-\int_{X}\langle A_\infty \nabla u_\infty, \nabla w \rangle dm_\infty \Bigr| \overset{k \to \infty}\to 0.
\end{align*} 
Next we show (II)$_k \to 0$ as $k \to \infty$. 
Combining $\sup_{n \in \N}|\b_i^n|<\infty$ with $L^2$-convergence of $\b_i^n$ to $\b_i^\infty$,  we have $\b^{n_k}_2 \to \b^{\infty}_2$ strongly in $L^2$, especially we have $\b^{n_k}_2(w) \to \b^{\infty}_2(w)$ strongly in $L^2$ . Since $u_k\to u_\infty$ weakly in $L^2$, we have that (II)$_k \to 0$. The quantity (IV)$_k \to 0$ in the same proof. 
Since $|{\rm div}\b_i^n|$ is uniformly bounded in $n$ and ${\rm div} \b_i^{n_k} \to {\rm div} \b_i^{\infty}$ in $L^2$, the quantity (III)$_k \to 0$ also goes to zero. 
 It is easy to check that (V)$_k \to 0$. The convergence of the dual forms can be shown {in} the same manner. 
\qed

\section{Convergence of Finite-dimensional Distributions} \label{sec: CFDDNS}
In this section, we show the weak convergence of finite-dimensional distributions. 
Recall that we identify $\iota_n(X_n)$ with $X_n$ and we omit $\iota_n$ for simplifying the notation. 
We first show the weak convergence of finite-dimensional distributions under the assumptions in Theorem \ref{thm: NSC2} in the case that the initial distribution is the Dirac measure $\delta_{\x_n}$.
 \begin{lem} \label{lem: FDC1}{\bf (Convergence of Finite-dimensional Distributions)}
Suppose the conditions assumed in Theorem \ref{thm: NSC2}. 
Then, for any $k \in \N$, $0 =t_0 < t_1 <t_2< \cdot \cdot \cdot <t_k<\infty$ and $f_1, f_2, ..., f_k \in C_b(X)$, the following holds:
\begin{align} \label{eq: FDC-RCD-4}
\mathbb E^{\x_n}[f_1(S^n_{t_1})\cdot\cdot\cdot f_k(S^n_{t_k})] \overset{n \to \infty}\to \mathbb E^{\x_\infty}[f_1(S^\infty_{t_1})\cdot\cdot\cdot f_k(S^\infty_{t_k})].
\end{align} 
For the dual process $\hat{\mathbb S}^{\x}_n$, the same statement holds.
\end{lem}
\proof
We omit the proof for the dual process which is the same as that of $\mathbb S^{\x}_n$. 
{Let $\{T_t\}_{t \ge 0}$ be the semigroup associated with the Dirichlet form $(\E, \F)$ corresponding to \eqref{form: non-sym}. According to the Gaussian heat kernel estimate \cite[Theorem 5.4]{LS17} (see also \eqref{GE}), under the assumptions in Theorem \ref{thm: NSC2}, we can easily show that $\{T^n_t\}_{t \ge 0}$ is a Feller semigroup, which implies the uniqueness of the corresponding diffusions for every starting point.}
\replaced{Therefore,}{Recall that} we have the following equality: for every $f \in C_b(X) \cap L^2(X;m_\infty)$,
\begin{align} \label{eq: KKT}
\mathbb E_n^{x}(f(S_t^n))=T^n_tf(x),
\end{align}
for {\it every} $x \in X_n$. 
\deleted{Note that by the heat kernel estimate \eqref{GE}, we can easily show that $\{T^n_t\}_{t \ge 0}$ is a Feller semigroup, which implies the uniqueness of the corresponding diffusions for every starting point.
Here $\{T^n_t\}_{t \ge 0}$ is the semigroup associated with $(\E_n, \F_n)$.}
By using the Markov property, for all $n \in \EN$, we have 
\begin{align*} 
&\mathbb E_n^{\x_n}[f_{1}(S^{n}_{t_1})\cdot\cdot\cdot f_k(S^{n}_{t_k})]  \notag
\\
&=T^n_{t_1-t_0}\Bigl(f_{1} T^n_{t_2-t_1} \Bigl(f_{2}\cdot \cdot \cdot T^n_{t_k-t_{k-1}}f_{k} \Bigr)\Bigr)(\x_n) \notag
\\
&=:\mathcal P_k^n(\x_n).
\end{align*}
By \cite[Corollary 4.18]{LS17}, the action of the semigroup $T_t^nf$ for $f \in L^\infty(m_n)$ is a H\"older continuous function whose H\"older constant and exponent are independent of $n$ (depending only on $N,K,D, \sup_{n \in \N}\||A_n|\|_\infty,  \sup_{n \in \N}\||\b_i^n|\|_\infty,  \sup_{n \in \N}\|c_n\|_{\infty}$).

For later arguments, we extend $\mathcal P_k^n$ to the whole space $X$ \added{by the McShane extension (\cite[Corollary 1,2]{Mc34})} (note that $\mathcal P_k^n$ is defined only on each $X_n$). 
The key point is to extend $\mathcal P_k^n$ to the whole space $X$ preserving its H\"older regularity and bounds.
Let $\widetilde{\mathcal P}_k^n$ be the following function on the whole space $X$:
\begin{align} \label{eq: EX}
\widetilde{\mathcal P_k^n}(x) :=\Bigl( \sup_{a \in X_n}\{\mathcal P_k^n(a)-H{d(a,x)^{\beta}}\} \wedge \sup_{a \in X_n} \mathcal P_k^n(a) \Bigr)  \vee  \inf_{a \in X_n} \mathcal P_k^n(a) \quad x \in X,
\end{align}
whereby $H$ and $\beta$ are the H\"older constant and exponent of the original function $\mathcal P_k^n$.
Then we have that $\widetilde{\mathcal P_k^n}$ is a bounded H\"older continuous function on the whole space $X$ with the same H\"older constant $H$, exponent $\beta$, and the same bound $\|\widetilde{\mathcal P_k^n}\|_\infty$, and satisfies $\tilde{\mathcal P_k^n}=\mathcal P_k^n$ on $X_n$.

Coming back to the proof of Lemma \ref{lem: FDC1}, we have that 
\begin{align}
&\Bigl| \mathbb E_n^{\x_n}[f_{1}(S^{n}_{t_1})\cdot\cdot\cdot f_k(S^{n}_{t_k})] -\mathbb E_n^{\x_\infty}[f_{1}(S^{\infty}_{t_1})\cdot\cdot\cdot f_k(S^{\infty}_{t_k})]\Bigr| \notag
\\
&=|\mathcal P_k^n(\x_n)-{\mathcal P}_k^\infty(\x_\infty)| \notag
\\
& \le |{\mathcal P}_k^n(\x_n)-\tilde{\mathcal P_k^n}(\x_\infty)|+|\tilde{\mathcal P_k^n}(\x_\infty)-\mathcal P_k^\infty(\x_\infty)| \notag
\\
&=:({\rm I})_n+({\rm II})_n.
\end{align}
Thus, it suffices to show that (I)$_n$ and (II)$_n$ converge to zero as $n$ goes to infinity. 
We first show that (I)$_n$ converges to zero as $n$ goes to infinity. 
By \replaced{the McShane extension}{Proposition \ref{prop: McS}}, we have 
\begin{align*}
({\rm I})_n= |{\mathcal P_k^n}(\x_n)-\tilde{\mathcal P_k^n}(\x_\infty)|
&= |\tilde{\mathcal P_k^n}(\x_n)-\tilde{\mathcal P_k^n}(\x_\infty)|
\\
&\le Hd(\x_n,\x_\infty)^{\beta}
\\
& \to 0 \quad (n \to \infty).
\end{align*}

Now we show that (II)$_n$ goes to zero as $n$ tends to infinity. Since 
$$\|T_t^nf\|_\infty =\|f\|_\infty \|\int_{X_n}p_n(t,x,y)m_n(dy)\|_\infty \le \|f\|_\infty,$$ for any $f \in C_b(X_n) \cap L^2(X;m_\infty)$, we have 
\begin{align} \label{ineq: BDDSG}
\sup_{n \in \N}\|\mathcal P_k^n\|_\infty \le \prod_{i=1}^k \|f_i\|_\infty<\infty.
\end{align}
Therefore, by the property of the McShane extension\deleted{ in Proposition \ref{prop: McS}}, we also have that 
\begin{align} \label{uniform bound of HS}
\sup_{n \in \N}\|\tilde{\mathcal P_k^n}\|_\infty<\infty.
\end{align}
By the uniform boundedness \eqref{uniform bound of HS} and the equi-continuity of $\{\tilde{\mathcal P_k^n}\}_{n \in \N}$, we can apply the Ascoli--Arzel\'a theorem to $\{\tilde{\mathcal P_k^n}\}_{n \in \N}$ so that $\{\tilde{\mathcal P_k^n}\}_{n \in \N}$  is relatively compact with respect to the uniform convergence so that for any subsequence $\{\tilde{\mathcal P_k^{n'}}\}_{\{n'\} \subset \{n\}}$, there exists a further subsequence $\{\tilde{\mathcal P_k^{n''}}\}_{\{n''\} \subset \{n'\}}$ satisfying 
\begin{align} \label{convergence: HS}
\tilde{\mathcal P_k^{n''}} \to F'' \quad \text{uniformly in}\ X.
\end{align}

On the other hand, we have that ${\mathcal P}_k^n$ converges to ${\mathcal P}_k^\infty$ $L^2$-strongly in the sense of Definition \ref{defn: Weak}. 
We give a proof below. 
\begin{lem} \label{lem: SCCHS}
${\mathcal P}_k^n$ converges to ${\mathcal P}_k^\infty$ in the $L^2$-strong sense in Definition \ref{defn: Weak}. 
\end{lem}
\proof
By Theorem \ref{thm: EQUIVHINO}, \ref{thm: HINOC}, the statement is true for $k=1$. Assume that the statement is true when $k=l$. By noting 
$$\mathcal P_{l+1}^n=T^n_{t_{l+1}-t_l}(f_{l+1}^{(n)}\mathcal P_l^n),$$
by Theorem \ref{thm: EQUIVHINO}, \ref{thm: HINOC} it suffices to show $f_{l+1}\mathcal P_l^n \to f_{l+1}\mathcal P_l^\infty$ strongly in $L^2$. This is easy to show because $\mathcal P_l^n \to \mathcal P_l^\infty$ strongly (the assumption of induction), $f_{l+1} \in C_b(X)$ and $\mathcal P_l^n$  is bounded uniformly in $n$ because of \eqref{uniform bound of HS}. Thus, the statement is true for any $k \in \N$. 
\qed

{\it Proof of Lemma \ref{lem: FDC1} {(Conclusion)}.}
By using Lemma \ref{lem: SCCHS} and the uniform convergence \eqref{convergence: HS}, it is easy to see that 
$$F''|_{X_\infty}=\mathcal P_k^\infty,$$
whereby $F''|_{X_\infty}$ means the restriction of $F''$ into $X_\infty$. 
The R.H.S. $\mathcal P_k^\infty$ of the above equality is clearly independent of choices of subsequences and thus, the limit $F''|_{X_\infty}$ is independent of choice of subsequences. Therefore, we conclude that 
\begin{align} \label{concl: UC}
\tilde{\mathcal P_k^n} \to \mathcal P_k^\infty \quad \text{uniformly in } \ X_\infty.
\end{align}
Going back to showing that (II)$_n$ goes to zero, we have that 
\begin{align*}
({\rm II})_n=|\tilde{\mathcal P_k^n}(\x_\infty)-\mathcal P_k^\infty(\x_\infty)|
&\le \|\tilde{\mathcal P_k^n}-\mathcal P_k^\infty\|_{\infty, X_\infty}
\\
& \to 0 \quad (n \to \infty).
\end{align*}
Here $\|\cdot\|_{\infty, X_\infty}$ means the uniform norm on $X_\infty$. Thus, we finish the proof of Lemma \ref{lem: FDC1}.
\qed

We now show the weak convergence of finite-dimensional distributions of $\mathbb S_n$ under Assumption \ref{asmp: NS3} with initial distributions $\nu_n$.
 Let us recall that $\zeta_{\mathbb S^n}$ and $\zeta_{\hat{\mathbb S}^n}$ denote lifetimes for $\mathbb S^n$ and $\hat{\mathbb S}^n$ respectively. Let $\zeta^n:=\min\{\zeta_{\mathbb S^n}, \zeta_{\hat{\mathbb S^n}}\}.$ 
\begin{lem} \label{lem: FDC-2}
Under Assumption \ref{asmp: NS3}, for any $k \in \N$, $0 =t_0 < t_1 <t_2< \cdot \cdot \cdot <t_k<\zeta_\infty<\infty$ and $f_1, f_2, ..., f_k \in C_b(X)$, the following holds:
\begin{align} \label{eq: FDC-RCD-4}
\mathbb E^{\nu_n}[f_1(S^n_{t_1})\cdot\cdot\cdot f_k(S^n_{t_k})] \overset{n \to \infty}\to \mathbb E^{\nu_\infty}[f_1(S^\infty_{t_1})\cdot\cdot\cdot f_k(S^\infty_{t_\infty})].
\end{align} 
For the dual process $\hat{\mathbb S}$, the same statement holds.  
\end{lem}
\proof
We omit the proof for $\hat{\mathbb S}$ since the proof is the same as the case of ${\mathbb S}$.
Since the limit diffusion $\mathbb S_\infty^{\x_\infty}$ is conservative,  it suffices to show the statement only for $f_1, f_2,...,f_k \in C_b(X) \cap L^2(X;m_\infty)$. In fact, for any $\e>0$ and $\zeta_\infty>T>0$, there exists $R=R(\e,T)$ so that the open ball $B_R(\x_\infty)$ satisfies
$$\mathbb P^{\nu_\infty}(S_t^\infty \in B_R(\x_\infty)^c) <\e \quad \forall t \in [0,T],$$
whereby $A^c:=X_\infty \setminus A$.  
By the strong $L^2$-convergence of the semigroup $\{T^n_t\}_{t \ge 0}$ following from Theorem \ref{thm: HINOC}, 
we have that 
\begin{align*}
 \lim_{n \to \infty}\mathbb P^{\nu_n}(S_t^n \in B_R(\x_\infty)^c) 
 < \e \quad \forall t \in [0,T].
\end{align*}
Therefore, for any $f_1,...,f_k \in C_b(X)$, for arbitrarily small $\delta>0$, we can take $R>0$ so that 
\begin{align*}
&\lim_{n \to \infty} \mathbb E^{\nu_n}(f_1(S_{t_1}^n)\cdot\cdot\cdot f_k(S_{t_k}^n))
\\
&= \lim_{n \to \infty} \mathbb E^{\nu_n}\Bigl(f_1(S_{t_1}^n)\cdot\cdot\cdot f_k(S_{t_k}^n): \ \bigcap_{j=1}^k\{S_{t_j}^n \in B_R(\x_\infty)\}\Bigr) 
\\
& \qquad \quad +  \lim_{n \to \infty} \mathbb E^{\nu_n}\Bigl(f_1(S_{t_1}^n)\cdot\cdot\cdot f_k(S_{t_k}^n): \ \Bigl(\bigcap_{j=1}^k\{S_{t_j}^n \in B_R(\x_\infty)\}\Bigr)^c\Bigr)
\\
&= \lim_{n \to \infty} \mathbb E^{\nu_n}\Bigl(f_1\1_{B_R}(S_{t_1}^n)\cdot\cdot\cdot f_k\1_{B_R}(S_{t_k}^n)\Bigr) + \delta.
\end{align*}
\added{Here we mean that, for an event $A \subset \Omega$, $\mathbb E^x(f(S_t): A):=\int_{\Omega \cap A}f(S_t(\omega))\mathbb P^{x}(d\omega)$.} 
Thus, we may show the proof only for  $f_1, f_2,...,f_k \in C_b(X) \cap L^2(X;m_\infty)$.
Since $\nu_n$ converges weakly to $\nu_\infty$ in $\mathcal P(X)$, for any $\e>0$, there exists a compact set $ K \subset X$ so that 
$$\sup_{n \in \N}\nu_n(K^c)<\e.$$
Thus, by \eqref{ineq: BDDSG}, for any $\delta>0$, there exists a compact set $ K \subset X$ so that 
\begin{align} \label{ineq: SFDC}
\sup_{n \in \N}\Bigl| \int_{X_n}\mathcal P_k^n d\m_n-\int_{K}\mathcal P_k^n\ d\nu_n\Bigr|& =  \sup_{n \in \N}\Bigl| \int_{X_n}\mathcal P_k^n (\1_{X_n}-\1_{K\cap X_n})d\nu_n \Bigl| \notag
\\
& \le \sup_{n \in \N}\|\mathcal P_k^n\|_{2,n}^{1/2} \nu_n(K^c) \notag
\\
&\le \Bigl(\prod_{i=1}^k\|f_i\|_{2,n}\Bigr) \sup_{n \in \N}\nu_n(K^c)< \delta.
\end{align}
We note that $\sup_{n\in \N}\prod_{i=1}^k\|f_i\|_{2,n}<\infty$ because $f_1, f_2,...,f_k \in C_b(X) \cap L^2(X;m_\infty)$ and $\mathcal X_n$ converges to $\mathcal X_\infty$ in the pmG sense. 
Take $r>0$ so that $K \subset B_r(\x_n):=\{x \in X: d(\x_n,x)<r\}$. Let $\tilde{\1}_r^R$ denote the following function: ($r<R$)
\begin{align*}
\tilde{\1}_r^R(x)=
\begin{cases} \dis
1 \quad \hspace{26mm} \text{$x \in B_r(\x_n)$},
\\
\dis 1-\frac{d(x, B_r(\x_n))}{R-r} \quad \text{$x \in B_R(\x_n) \setminus B_r(\x_n)$},
\\
0 \hspace{30mm} \text{o.w.}
\end{cases}
\end{align*}
Then $\tilde{\1}_r^R \in C_{bs}(X)$.
Thus, by Theorem \ref{thm: EQUIVHINO}, \ref{thm: HINOC} and \eqref{ineq: SFDC}, for any $\delta>0$, there exists $r>0$ so that
\begin{align*}
& \Bigl| \mathbb E^{\nu_n}[f_1(S^n_{t_1})\cdot\cdot\cdot f_k(S^n_{t_k})]  - \mathbb E^{\nu_\infty}[f_1(S^\infty_{t_1})\cdot\cdot\cdot f_k(S^\infty_{t_\infty})] \Bigr|
\\
&=  \Bigl| \int_{X_n}\mathcal P_k^n d\nu_n - \int_{X_\infty}\mathcal P_k^\infty d\nu_\infty \Bigr|
\\
&=  \Bigl| \int_{X_n}\mathcal P_k^n d\nu_n - \int_{X_n}\tilde{\1}_r^R \mathcal P_k^n d\nu_n+ \int_{X_n}\tilde{\1}_r^R\mathcal P_k^n \phi_n^1dm_n- \int_{X_n}\tilde{\1}_r^R\mathcal P_k^\infty \phi_\infty^1dm_\infty
\\
& \quad +\int_{X_n} \tilde{\1}_r^R \mathcal P_k^\infty d\nu_\infty- \int_{X_\infty}\mathcal P_k^\infty d\nu_\infty \Bigr|
\\
& \le \delta+ \Bigl|\int_{X}\tilde{\1}_r^R \mathcal P_k^n \phi_ndm_n- \int_{X}\tilde{\1}_r^R  \mathcal P_k^\infty \phi_\infty dm_\infty\Bigr| + \delta
\\
& \overset{n \to \infty}\to 2\delta. 
\end{align*} 
Here, in the fifth line above, in the first $\delta$,  we used \eqref{ineq: SFDC} and in the second $\delta$, we used the tightness of the single measure $m_\infty$. The middle term in the fifth line converges to zero because of the $L^2$-strong convergence of the semigroup $T^n_t$ by Theorem \ref{thm: EQUIVHINO}  and Theorem \ref{thm: HINOC}, and $L^2$-weak convergence of $\phi_n$ to $\phi_\infty$. 
Thus, we have completed the proof.
\qed

\section{Tightness} \label{sec: TGHTNS}
In this section, we investigate the tightness of the diffusion processes associated with $(\E_n, \F_n)$.
According to  \cite{T89} and \cite{T08}, we have a decomposition of additive functionals for non-symmetric forms, which is called {\it Lyons-Zheng decomposition}. Suppose Assumption \ref{asmp: NS1}. Let $\mathbb S=(\Omega, \{\mathcal M_t\}_{t \ge 0}, \{S_t\}_{t \ge 0}, \{\mathbb P^x\}_{x \in X})$ be a diffusion process on $\Omega=D([0,\infty);X_{\partial})$ associated with the Dirichlet form $(\E, \F)$ corresponding to \eqref{form: non-sym}. We take $S_t(\omega)=\omega(t)$ as a coordinate process for $\omega \in \Omega$. Let $\hat{\mathbb S}=(\Omega, \{\mathcal M_t\}_{t \ge 0}, \{S_t\}_{t \ge 0}, \{\hat{\mathbb P}^x\}_{x \in X})$ be a dual process associated with the dual form $(\hat{\E}, \F)$. Let $r_T$ be a time reversal operator defined as follows:
\begin{align*}
r_s(\omega)(t)=
\begin{cases}
\omega((s-t)-), \quad 0 \le t \le s<\infty,
\\
\omega(0), \hspace{16mm} s <t.
\end{cases}
\end{align*}
Here $\omega(t-):=\lim_{s \uparrow t}\omega(s)$. Since $(\E, \F)$ is local, the corresponding processes are diffusive and jump only to the cemetery point $\partial$. Thus, we may omit to write $\omega(t-)$ before lifetime and simply write $\omega(t)$.
By the Fukushima decomposition, we have that for $f \in \F$
\begin{align} \label{eq: FD1}
f(S_t)-f(S_0)=M_t^{[f]}-N_t^{[f]}, \quad \text{$\mathbb P^x$-a.e., q.e. $x$},
\end{align}
whereby $M^{[f]}$ is a martingale and $N^{[f]}$ is a zero-energy process. 
{Here we mean by {\it zero-energy} $e(N^{[f]})=0$, in which the energy of $N^{[f]}$ (generally, the energy of additive functionals) is defined as follows:
$$
e(N^{[f]}):=\lim_{\alpha \to \infty} \alpha^2 \EE^{m}\Bigl[\int_0^\infty e^{- \alpha t}(N_t^{[f]})^2dt\Bigr].
$$
}
Similarly,  we have that for $f \in \F$,
\begin{align} \label{eq: FD2}
f(S_t)-f(S_0)=\hat{M}_t^{[f]}-\hat{N}_t^{[f]}, \quad \text{$\hat{\mathbb P}^x$-a.e., q.e. $x$}.
\end{align}
Let $\zeta_{\mathbb S}$ and $\zeta_{\hat{\mathbb S}}$ be lifetimes for $\mathbb S$ and $\hat{\mathbb S}$ respectively. Let $\zeta:=\min\{\zeta_{\mathbb S}, \zeta_{\hat{\mathbb S}}\}.$ We note that, on $\{\zeta>T\},$ we have that, for an $\mathcal M_T$-measurable function $F$, 
$$\hat{\mathbb E}^m(F(r_T\omega))=\mathbb E^m(F(\omega)).$$
By 
\cite{T08}, for $f \in \mathcal F$, we have that on $\{\zeta > T\}$,
\begin{align} \label{eq: tight2}
\tilde{f}(S_t)-\tilde{f}(S_0)&=\frac{1}{2}M_t^{[f]}-\frac{1}{2}(\hat M_T^{[f]}(r_T)-\hat M_{T-t}^{[f]}(r_T)) +\frac{1}{2}(N_t^{[f]}-\hat{N}_t^{[f]}),
\end{align}
for $0 \le t \le T$ $\mathbb P^{m}$-a.e. Here $\tilde{f}$ means a quasi-continuous modification of $f$.

Now we estimate $\frac{1}{2}(N_t^{[u]}-\hat{N}_t^{[u]})$.
\begin{lem} \label{lem: NSMT1}
\added{Suppose Assumption \ref{asmp: NS1}, $|A|, |b_1|, |b_2|, c \in L^\infty(X,m)$, $\b_i \in {\rm Div}^2(X,d,m)$ ($i=1,2$) and symmetry of $A$.}
For $f \in \mathcal F$ and $t \ge 0$, it holds that on $\{\zeta >T\}$,
$$\hat N_t^{[f]}-N_t^{[f]}=\int_0^t \Bigl(2\b_1(f)-2\b_2(f) -f{\rm div} \b_1 + f{\rm div} \b_2\Bigr)(S_s)ds,$$
for $0 \le t \le T$ $\mathbb P^{m}$-a.e.
\end{lem}
\proof
First we prove the statement for $f  \in \mathcal D(\hat{{L}})$, whereby $\hat{{L}}$ denotes the generator associated with $(\hat \E, \F)$ and $\mathcal D (\hat{{L}})$ denotes the domain of $\hat{{L}}$.
In this case, we have $\hat N_t^{[f]}=\int_0^t \hat{{L}}f(X_s)ds$ and thus, we see $\hat N_t^{[f]}(r_t)=\hat N_t^{[f]}.$
Then for $f \in \mathcal F$, we have on $\{\zeta >T\}$,
\begin{align*}
\EE^{gm}[\hat N_t^{[f]}]&=\hat  \EE^m[\hat N_t^{[f]}(r_t)\tilde g(S_t)]
\\
&= \hat \EE^m[\hat N_t^{[f]}\tilde g(S_t)]
\\
&=\hat \EE^m[\hat N_t^{[f]}\tilde g(S_t)]+\hat \EE^m[\hat N_t^{[f]}(\tilde g(S_t)-\tilde g(S_0))], \quad \forall g \in \mathcal F,
\end{align*}
whereby $\tilde{g}$ denotes a quasi-continuous modification of $g$.
\replaced{We have}{Since} 
\begin{align*}
&\alpha^2\Bigl|\hat \EE^m\Bigl[\int_0^\infty e^{-\alpha t}\hat N_t^{[f]}(\tilde g(S_t)-\tilde g(S_0))\Bigr]\Bigr|
\\
& \le \Bigl(\alpha^2 \hat \EE^m\Bigl[\int_0^\infty e^{-\alpha t}(\hat N_t^{[f]})^2dt\Bigr]\Bigr)^{1/2}\Bigl(\alpha^2 \hat \EE^m\Bigl[\int_0^\infty e^{-\alpha t}(\tilde g(X_t)-\tilde g(X_0))^2dt\Bigr]\Bigr)^{1/2} 
\\
& \overset{\alpha \to \infty}\to \hat e(\hat N^{[f]})^{{1/2}}\hat e({\tilde g(X_t)-\tilde g(X_0)})^{1/2}=0.
\end{align*}
\added{Therefore,} by \cite[Theorem 5.3.1]{O13}, we have
\begin{align*}
&\lim_{ \alpha \to \infty}\alpha^2\EE^{gm}\Bigl[\int_0^\infty e^{- \alpha t}\hat N_t^{[f]}dt\Bigr]
\\
&=\lim_{ \alpha \to \infty}\alpha^2\EE^{m}\Bigl[\int_0^\infty e^{- \alpha t}\hat N_t^{[f]}\tilde{g}(S_t)dt\Bigr]
\\
&=\lim_{ \alpha \to \infty}\alpha^2\hat{\EE}^{m}\Bigl[\int_0^\infty e^{- \alpha t}\hat N_t^{[f]}\tilde{g}(S_t)dt\Bigr]+\alpha^2\Bigl|\hat \EE^m\Bigl[\int_0^\infty e^{-\alpha t}\hat N_t^{[f]}(\tilde g(S_t)-\tilde g(S_0))\Bigr]\Bigr|
\\
&=\lim_{ \alpha \to \infty}\alpha^2 \hat \EE^{gm}\Bigl[\int_0^\infty e^{-\alpha t}\hat N_t^{f}dt\Bigr]
\\
&=-\hat \E(f,g).
\end{align*}
Since it holds that (recall symmetry of $A$) 
$$\hat \E(f,g)+\int_X\bigl(2\b_1(f)-2\b_2(f) -  f{\rm div} \b_1 + f{\rm div} \b_2\bigr)g dm=\E(f,g), $$
by \cite[Theorem 5.3.1]{O13}, we obtain 
\begin{align} \label{eq: rel. dual}
&\lim_{\alpha \to \infty} \alpha^2 \EE^{gm}\Bigl[\int_0^\infty \Bigl(e^{- \alpha t}\hat N_t^{[f]}-\int_0^t\Bigl(2\b_1(f)-2\b_2(f) -f{\rm div} \b_1 + f{\rm div} \b_2\Bigr)(S_s)ds\Bigr)dt\Bigr] \notag
\\
&=\lim_{\alpha \to 
\infty}\alpha^2 \EE^{gm}\Bigl[\int_0^\infty e^{- \alpha t}N_t^{[f]}dt\Bigr].
\end{align}
On the other hand, we can calculate the energy of $\hat{N}^{[f]}$ as follows:
\begin{align} \label{eq: energy-zero}
e(\hat N^{[f]})&=\lim_{\alpha \to \infty} \alpha^2 \EE^{m}\Bigl[\int_0^\infty e^{- \alpha t}(\hat N_t^{[f]})^2dt\Bigr] \notag
\\
&=\lim_{\alpha \to \infty} \alpha^2 \hat \EE^{m}\Bigl[\int_0^\infty e^{- \alpha t}(\hat N_t^{[f]})^2dt\Bigr] \notag
\\
&=\hat e(\hat N^{[f]}) \notag
\\
&=0.
\end{align}
Thus, $\hat N_t^{[f]}-\int_0^t (2\b_1(f)-2\b_2(f)-f{\rm div} \b_1 + f{\rm div} \b_2)(S_s) ds$ is an \replaced{additive functional}{AF} of $\mathbb S$ with zero energy. By \eqref{eq: rel. dual}, \eqref{eq: energy-zero}, and \cite[Theorem 5.3.1]{O13},  we have the desired result for $f \in \mathcal D(\hat A).$

For general $f \in \mathcal F$, we can take a sequence $f_n \in \mathcal D(\hat{{L}})$ so that $f_n$ converges to $f$ with respect to the norm of the symmetric part $\tilde{\mathcal E}_{1}$ and for  q.e.\ $x$ (\cite[Theorem 5.1.3]{O13}), 
$$\hat{\mathbb P}^x[\hat{\Gamma}_T]=1,$$
where 
$$\hat \Gamma_T=\{\omega \in \Omega: \hat N_t^{[f_n]}(\omega) \ \text{converges to } \hat N_t^{[f]}(\omega) \text{\ uniformly in $t$ on $[0,T]$}\}.$$
Since we have that on $\{\zeta >T\}$,  
\begin{align} \label{formula: ZE}
\hat N_{t}^{[f_n]}(r_T \omega)&= \int_{T-t}^T\hat{{L}}f_n(S_s(\omega))ds =\hat N_{T}^{[f_n]}(\omega)-\hat N_{T-t}^{[f_n]}(\omega), 
\end{align}
the set $\hat \Gamma_T$ is $r_T$-invariant, i.e., $\{r_T \omega \in \hat \Gamma_T\}= \hat \Gamma_T$. Therefore, the complement $\hat \Gamma_T^c$ of $\hat \Gamma_T$ is also $r_T$-invariant. Thus, we obtain 
\begin{align*}
\mathbb P^m(\hat{\Gamma}_T^c)=\hat{\mathbb P}^m(r_T \omega \in \hat{\Gamma}_T^c)
=\hat{\mathbb P}^m(\hat{\Gamma}_T^c)=0.
\end{align*}
Therefore, we can conclude the desired result for general $f \in \mathcal F$\deleted{ by uniform approximation}. 
\qed

\added{By the previous lemma, by easy calculation, we have that,}
for $f \in \mathcal F$, on $\{\zeta>T\}$, 
\begin{align} \label{eq: LZR}
\tilde{f}(S_s)-\tilde{f}(S_0)&=\frac{1}{2}M_t^{[f]}-\frac{1}{2}(\hat M_T^{[f]}-\hat M_{T-t}^{[f]}(r_T)) \notag
\\
&\quad -\int_0^t \Bigl(\b_1(f)-\b_2(f)-\frac{1}{2}f{\rm div} \b_1+ \frac{1}{2}f{\rm div} \b_2\Bigr)(S_s) ds,
\end{align}
for $0 \le t \le T$ $\mathbb P^m$-a.e.

\begin{lem} \label{thm: tightness-nonsym}
Under Assumption \ref{asmp: NS3}, 
 $\{\mathbb S^{\nu_n}\1_{\{\zeta^n >T \}}\}_{n \in \N}$ and  $\{\hat{\mathbb S}^{\nu_n}\1_{\{\zeta^n >T \}}\}_{n \in \N}$ are tight in $\mathcal P_{\le 1}(C([0,T]; X))$ for any $T>0$.
\end{lem}
\proof
We only show the tightness of $\{\mathbb S^{\nu_n}\}_{n \in \N}$ since the proof for the dual processes is the same. 
Let us denote the law of $h(B^n)$ for $h \in {\rm Lip}_{bs}(X)$ as follows:
$$
 \mathbb S^{\nu_n, h}=(h(S^n), \mathbb P_n^{\nu_n}).
$$
Here we set $h(\partial)=0.$
It is easy to show that ${\rm Lip}_{bs}(X)$ strongly separates points in $C_b(X)$, that is, for every $x$ and $\e>0$, there exists a finite set $\{h_i\}_{i=1}^l \subset {\rm Lip}_{bs}(X)$ so that 
$$\inf_{y: d(y,x) \ge \e}\max_{1 \le i \le l}|h_i(x)-h_i(y)|>0.$$
By \cite[Theorem 3.9.1, Corollary  3.9.2]{EK86}  (we can apply these statements also to the space $\mathcal P_{\le 1}(C([0,T]; X))$ of sub-probability measures) and Lemma \ref{lem: FDC-2}, the following (i) follows from (ii): For any $T>0$,
\begin{enumerate}\label{conv: EK}
\item[(i)]  $\{\mathbb S^{\nu_n}\1_{\{\zeta^n >T \}}\}_{n \in \N}$ is tight in $\mathcal P_{\le 1}(C([0,T]; X))$;
\item[(ii)]  $\{\mathbb S^{\nu_n,h}\1_{\{\zeta^n >T \}}\}_{n \in \N}$ is tight in $\mathcal P_{\le 1}(C([0,T]; \R))$ for $\forall h \in {\rm Lip}_{bs}(X)$.
\end{enumerate}
In fact, we can show the compact containment condition \cite[(9.1) in Theorem 3.9.1]{EK86} according to (ii) and Lemma \ref{lem: FDC-2} by a proof similar to \cite[Corollary  3.9.2]{EK86}.
We note that, although \cite[Theorem 3.9.1, Corollary 3.9.2]{EK86} gives sufficient conditions for tightness only in the c\`adl\`ag space $D([0,T];X)$, since the laws of each diffusions $\mathbb S_n^{\nu_n}$ and $\mathbb S_\infty^{\nu_\infty}$ have their support on the space of continuous paths $C([0,T];X)$ before lifetime because of the locality of $(\E_n, \F_n)$, the tightness in $D([0,T];X)$ implies the tightness in $C([0,T]; X)$. See, e.g., \cite[Lemma 5 in Appendix]{FK97} for this point.
Thus, we will show that (ii) holds, i.e., for any $T$ and  any $h \in {\rm Lip}_{bs}(X)$, 
 \begin{align}  \label{conv: EK}
\{\mathbb S^{\nu_n,h}\1_{\{\zeta^n >T \}}\}_{n \in \N} \quad \text{ is tight in} \quad \mathcal P_{\le 1}(C([0,T];X)).
 \end{align}
Since $\nu_n$ converges weakly to $\nu_\infty$ in $\mathcal P(X)$, 
the laws of the initial distributions $\{{(h(S^n_0)}, \mathbb P_n^{\nu_n})\}_{n \in \N}=\{h_\#\nu_n\}_{n \in \N}$ {are} clearly tight in $\mathcal P(\R)$.
For $\delta>0$, let us define
 $$
 L_{\eta,T}^{n,h}(x):=\mathbb P_n^x\Bigl(\sup_{\substack{0 \le s,t\le T \\ |t-s|\le \eta}}|h(S^n_t)-h(S^n_s)| >\delta: \{\zeta^n >T\}\Bigr).
 $$
By the local property of $(\E_n, \F_n)$, we see that $S^n$ is continuous in the event $\{\zeta^n >T\}$. Thus, the desired result we would like to show is the following:
 \begin{align} \label{convergence: tightness-RCD}
 \lim_{\eta\to 0}\sup_{n \in \N}\int_{X_n}L_{\eta,T}^{n,h}\ \nu_n(dx)=0,
 \end{align}
 for any $T>0$.
For any $\e>0$, there exists $R>0$ so that 
 \begin{align*}
 \int_{X_n}L_{\eta,T}^{n,h}\nu_n(dx)
 &=\|\phi_n\1_{B_R(\x_n)}\|_\infty \int_{X_n}L_{\eta,T}^{n,h}\1_{B_R(\x_n)} dm_n+\nu_n(B_R^c(\x_n))
 \\
 &<\|\phi_n\1_{B_R(\x_n)}\|_\infty \int_{X_n}L_{\eta,T}^{n,h}\1_{B_R(\x_n)} dm_n+\e.
 \end{align*}
Let $m_{n,R}:=\1_{B_{R}(\x_n)}m_n$. 
We have 

\begin{align} \label{eq: tight1}%
\int_{X_n}L_{\eta,T}^{n,h}\ dm_{n,R}  & =\mathbb P_{n}^{m_{n,R}}\Bigl(\sup_{\substack{0 \le s,t\le T \\ |t-s|\le h}}|h(S^n_t)-h(S^n_s)| >\delta: \{\zeta^n >T\} \Bigr)
\\
&:= {\rm (I)}_{{n},\eta}. \notag
\end{align}
It suffices to show that, for any $T,R>0$, 
$$\sup_{n \in \N}{\rm (I)}_{{n},\eta} \to 0, \quad \eta \to 0.$$
By \replaced{the equality \eqref{eq: LZR}}{Proposition \ref{prop: NSLZ1}}, we have that, on $\{\zeta^n>T\}$, 
\begin{align} \label{eq: tight2-1}
h(S^n_s)-h(S^n_0)&=\frac{1}{2}M_t^{[h],n}-\frac{1}{2}(\hat M_T^{[h],n}-\hat M_{T-t}^{[h],n}(r_T)) \notag
\\
&\quad -\int_0^t \Bigl(\b^n_1(h)-\b^n_2(h)-\frac{1}{2}h{\rm div} \b^n_1+ \frac{1}{2}h{\rm div} \b^n_2\Bigr)(S^n_s) ds,
\end{align}
for $0 \le t \le T$ $\mathbb P^m$-a.e.
Thus, we have 
\begin{align}\label{eq: tight3-1} 
{\rm (I)}_{n,\eta}&=\mathbb P^{m_{n,R}}\Bigl(\sup_{\substack{0 \le s,t\le T \\ |t-s|\le \eta}}|h(S^n_t)-h(S^n_s)| >\delta: \{\zeta^n >T\}\Bigr) 
\notag
\\
&\le \mathbb P^{m_{n,R}}\Bigl(\sup_{\substack{0 \le s,t\le T \\ |t-s|\le \eta}}\bigl|  M_t^{[h],n}- M_s^{[h],n} \bigr| > \delta: \{\zeta^n >T\}\Bigr) \notag
\\
& \quad + \mathbb P^{m_{n,R}}\Bigl(\sup_{\substack{0 \le s,t\le T \\ |t-s|\le \eta}}\bigl| \hat M_{T-t}^{[h],n}(r_T)-\hat M_{T-s}^{[h],n}(r_T) \bigr| > \delta: \{\zeta^n >T\}\Bigr) \notag
\\
&\quad +  \mathbb P^{m_{n,R}}\Bigl(\sup_{\substack{0 \le s,t\le T \\ |t-s|\le \eta}}\Bigl| \int_s^t \Bigl(\b^n_2(h)-\b^n_1(h)-\frac{1}{2}h{\rm div} \b^n_2 + \frac{1}{2}h{\rm div} \b^n_1\Bigr)(S_l)dl \Bigr| > \delta: \{\zeta^n >T\} \Bigr).
\end{align}
Noting that $\hat{\mathbb E}^m(F(r_T\omega))=\mathbb E^m(F(\omega))$ for an $\mathcal M_T$-measurable function $F$ on $\{\zeta^n>T\}$, we have 
\begin{align}\label{ineq: tight3}
&(\text{R.H.S. of \eqref{eq: tight3-1}}) = \mathbb P^{m_{n,R}}(\sup_{\substack{0 \le s,t\le T \\ |t-s|\le \eta}} \bigr| M_t^{[h],n}-M_s^{[h],n} \bigr|  > \delta: \{\zeta^n >T\}\Bigr) \notag
\\
& \quad + \hat{\mathbb P}^{m_{n,R}}\Bigl(\sup_{\substack{0 \le s,t\le T \\ |t-s|\le \eta}}\bigl| \hat M_{T-t}^{[h],n}-\hat M_{T-s}^{[h],n} \bigr| > \delta: \{\zeta^n >T\}\Bigr) \notag
\\
& \quad +  \mathbb P^{m_{n,R}}\Bigl(\sup_{\substack{0 \le s,t\le T \\ |t-s|\le \eta}}\Bigl| \int_s^t \Bigl(\b^n_2(h)-\b^n_1(h)-\frac{1}{2}h{\rm div} \b^n_2 + \frac{1}{2}h{\rm div} \b^n_1\Bigr)(S_l)dl \Bigr| > \delta: \{\zeta^n >T\} \Bigr)  \notag
 \\
& {\le} \mathbb P^{m_{n,R}}\Bigl(\sup_{\substack{0 \le s,t\le T \\ |t-s|\le \eta}} \bigr| M_t^{[h],n}-M_s^{[h],n} \bigr|  > \delta: \{\zeta^n >T\}\Bigr) \notag
 \\
& \quad + \hat{\mathbb P}^{m_{n,R}}\Bigl(\sup_{\substack{0 \le s,t\le T \\ |t-s|\le \eta}}\bigl| \hat M_{T-t}^{[h],n}-\hat M_{T-s}^{[h],n} \bigr| > \delta: \{\zeta^n >T\}\Bigr) \notag
\\
& \quad +  \mathbb P^{m_{n,R}}\Bigl((\||\b^n_1|\|_\infty+\||\b^n_2|\|_\infty+\|{\rm div} \b^n_1\|_\infty+\|{\rm div} \b^n_2\|_\infty)\eta  > \delta: \{\zeta^n >T\}\Bigr).
\end{align}
We first estimate the martingale part. 
Since $M^{[h],n}$ is a continuous martingale, by the martingale representation theorem, there exists \replaced{a}{the} one-dimensional Brownian motion  ${\mathbf B}^n(t)$ on an extended probability space $(\tilde{\Omega}, \tilde{\mathcal M}, \tilde{\mathbb P}_n^x)$, whereby $M^{[h],n}$ is represented as a time-changed Brownian motion with respect to the quadratic variation $\tilde{\mathbb P}_n^x$-a.s, q.e.\ $x \in X_n$ (see, e.g.,  Ikeda--Watanabe \cite[Chapter II Theorem 7.3']{IW89}). That is, for q.e.\ $x \in X_n$,
\begin{align} 
M^{[h],n}_t={\mathbf B}^n(\la M^{[h],n} \ra_t)={\mathbf B}^n\Bigl(\int_0^t \frac{d\mu^n_{\la h \ra}}{dm_n}(S_u^n)du\Bigr)={\mathbf B}^n\Bigl(\int_0^t \langle A_n\nabla h, \nabla h \rangle(S_u^n)du\Bigr)  \quad \text{$\tilde{\mathbb P}_n^x$-a.s.}
\end{align}
Here $\mu^n_{\la \cdot \ra}$ means the energy measure associated with ${\sf Ch}_n$: ${\sf Ch}_n(f)=\int_{X_n}\mu^n_{\la f \ra}(dx)$ for $f \in \mathcal F^{{n}}$. 
Since $|\nabla h| \le {\rm Lip}(h)$, we have
\begin{align}
&\{\omega \in \tilde{\Omega}:\sup_{\substack{0 \le s,t\le T \\ |t-s|\le \eta}}\bigr| M_t^{[h],n}-M_s^{[h],n} \bigr|  > \delta \} \notag
\\
&\subset\{\omega \in \tilde{\Omega}:\sup_{\substack{0 \le s,t\le T \\ |t-s|\le \eta}}\Bigr| {\mathbf B}^n\Bigl(\int_0^t  \||A_n|\|_\infty|\nabla h|^2(S_u^n)du\Bigr)-{\mathbf B}^n\Bigl(\int_0^s  \||A_n|\|_\infty|\nabla h|^2(S_u^n)du\Bigr) \Bigr|  > \delta \} \notag
\\
&\subset \{\omega \in \tilde{\Omega}:\sup_{\substack{0 \le s,t\le  \||A_n|\|_\infty{\rm Lip}(h)^2T \\ |t-s|\le \||A_n|\|_\infty{\rm Lip}(h)^2\eta}}\bigr| {\mathbf B}^n(t)-{\mathbf B}^n(s) \bigr|  > \delta \}. \notag
\end{align}
Let $\mathbb W$ be the standard Wiener measure on $C([0,\infty);\R)$. Let 
$$\theta(\eta,h):=\mathbb W_n(\sup_{\substack{0 \le s,t\le \||A_n|\|_\infty{\rm Lip}(h)^2T \\ |t-s|\le \||A_n|\|_\infty{\rm Lip}(h)^2\eta}}|\omega(t)-\omega(s)|>\delta).$$
By \eqref{ineq: tight3} and noting $\sup_{n \in \N}m_n(B_R(\x_n))<\infty$ because of the weak convergence of $m_n$, we have, for any $T>0$,
\begin{align} \label{tight: est1}
& \sup_{n \in \N}\mathbb P^{m_{n,R}}(\sup_{\substack{0 \le s,t\le T \\ |t-s|\le \eta}} \bigr| M_t^{[h],n}-M_s^{[h],n} \bigr|  > \delta: \{\zeta^n >T\}) \notag
 \\
 &\le \theta(\eta,h)\sup_{n \in \N}m_{n}(B_R(\x_n)) \notag
 \\
 & \overset{\eta \to 0}\to 0.
 \end{align}
 The dual martingale part can be estimated in a similar way, so, we omit the proof. 
For the remaining part, we can see that the following uniform estimate in $n$: 
$$\sup_{n \in \N}\mathbb P^{m_{n,R}}\Bigl((\||\b^n_1|\|_\infty+\||\b^n_2|\|_\infty+\|{\rm div} \b^n_1\|_\infty+\|{\rm div} \b^n_2\|_\infty)\eta  > \delta: \{\zeta^n >T\}\Bigr) {=}0,$$
{provided that 
$$\eta\sup_{n \in \N}(\||\b^n_1|\|_\infty+\||\b^n_2|\|_\infty+\|{\rm div} \b^n_1\|_\infty+\|{\rm div} \b^n_2\|_\infty)<\delta.$$
}
Thus, we obtain $\sup_{n \in \N}{\rm (I)}_{{n},\eta} \to 0$ as $\eta \to 0$. We have finished the proof. 
\qed

\begin{lem} \label{thm: tightness-nonsym1}
Suppose the conditions assumed in Theorem \ref{thm: NSC2}. 
Then 
 $\{\mathbb S^{\x_n}\}_{n \in \N}$ is tight in $\mathcal P(C([0,\infty), X))$.
\end{lem}
 \proof
Since $\x_n$ converges to $\x_\infty$ in $(X,d)$, the laws of the initial distributions $\{S^n_0\}_{n \in \N}=\{\delta_{x_n}\}_{n \in \N}$ \replaced{are}{is} clearly tight in $\mathcal P(X)$. 
Thus, it suffices to show the following (see \cite[Theorem 12.3]{B99}): 
for each $T>0$, there exist $\beta>0$, $C>0$ and $\theta>1$ such that, for all $n \in \N$,
	\begin{align} \label{criterion: tight5}
	&\mathbb E^{x_n}[\tilde{d}^{\beta}(S^{n}_t, S^{n}_{t+h})]
	\le Ch^{\theta}, \quad (\text{{for every }} 0 \le t \le T \ \text{and} \ 0 \le h \le 1), 
	\end{align}
	whereby $\tilde{d}(x,y):=d(x,y) \wedge 1$.
By the Markov property, we have 
\begin{align} \label{eq: tight5}
	&\text{L.H.S. of \eqref{criterion: tight5}}  \notag
	\\
	& =\int_{X_n \times X_n} p_n(t,x_n,y)p_n(h,y,z)\tilde{d}^\beta(\iota_n(y),\iota_n(z)) m_n(dy)m_n(dz). \notag
	\\
	& \le \int_{X_n \times X_n} p_n(t,x_n,y)p_n(h,y,z){d}^\beta(\iota_n(y),\iota_n(z))m_n(dy)m_n(dz).
\end{align}
By the generalized Bishop--Gromov inequality \cite[Proposition 3.9]{EKS15}, we have the following volume growth estimate: there exist positive constants $\nu=\nu(N,K, D)>0$ and $c=c(N,K,D)>0$ such that, for all $n \in \N$
\begin{align} \label{VE}
	m_n(B_r(x)) \ge cr^{2\nu} \quad  (0 \le  r \le 1 \wedge D).
\end{align} 

On the other hand, the volume doubling property (\cite{Sturm06}) and the Poincar\'e inequality (\cite{HK95, Raj12a, Raj12b}) hold under RCD$^*(K,N)$ condition. According to \cite[Theorem 5.4]{LS17} and \eqref{VE}, we have that
there exist positive constants $C_1$ and $C_2$ depending only on $N,K,D, \sup_{n \in \N}\{\||A_n|\|_\infty+\||\b_1^n|\|_\infty+\||\b_2^n|\|_\infty+\|c_n\|_\infty\}$ and $T$ so that 
 \begin{align} \label{GE}
p(t,x,y) \le \frac{C_1}{t^\nu}\exp\Bigl\{-C_2\frac{d(x,y)^2}{t}\Bigr\},
 \end{align}
for all $x, y \in X$ and $0<t\le D^{2}$. 
Thus, we have
\begin{align} 
	\int_{X_n} &p_n(s,y,z){d}^\beta(\iota_n(y),\iota_n(z)) m_n(dz) \notag
	\\
	&\le \frac{C_1}{cs^{\nu}}\int_{X_n}\exp\Bigl(- C_2\frac{d_n(y,z)^2}{s}\Bigr){d}^\beta(\iota_n(y),\iota_n(z))m_n(dz) \notag
	\\
	&\le \frac{C_1}{cs^{\nu}}\int_{X_n}\exp\Bigl(- C_2\frac{d_n(y,z)^2}{s}\Bigr){d}^\beta_n(y,z)m_n(dz) \notag
	\\
	& \le C_1c^{-1}C_2^{2/\beta}C_3s^{\beta/2-\nu}m_n(X_n)\sup_{y,z\in X_n}\Bigl\{\Bigl(C_2\frac{d_n(y,z)^2}{s}\Bigr)^{\beta/2}\exp\Bigl(- C_2\frac{d_n(y,z)^2}{s}\Bigr)\Bigr\} \notag
	\\ \label{ineq: UGEP}
	& \le  C_1c^{-1}C_2^{2/\beta}M_\beta s^{\beta/2-\nu}\notag
	\\
	& = C_4s^{\beta/2-\nu}.
\end{align}

By \eqref{ineq: UGEP}, we have
\begin{align}
	\text{R.H.S. of \eqref{eq: tight5}} & \le C_4h^{\beta/2-\nu} \int_{X_n} p_n(t,x_n,y)m_n(dy) \notag
	\\
	& \le C_4h^{\beta/2-\nu}.
\end{align}
Thus, we finish the proof by taking $\beta>0$ such that $\beta/2-\nu>1$, and set $\theta=\beta/2-\nu$.
\qed

\deleted{Now we prove Theorem \ref{thm: NSC1}, \ref{thm: NSC2}.}
By Lemma \ref{lem: FDC-2}, \ref{thm: tightness-nonsym}, we can finish the proof of Theorem \ref{thm: NSC1}. By Lemma \ref{lem: FDC1}, \ref{thm: tightness-nonsym1}, we finish the proof of Theorem \ref{thm: NSC2}.

\section{Conservativeness} \label{sec: CONSVNS}
In this section, under Assumption \ref{asmp: NS1} with $\alpha_0=0$, we give a criterion for the conservativeness of $(\E, \F)$ and $(\hat{\E}, \F)$.
In the case of finite mass $m(X)<\infty$, if
{
\begin{align*}
 {\rm div} \b_i=c, \quad i=1,2,
\end{align*}
}
then it is easy to check the conservativeness since $\1 \in \mathcal F$ and $\E(\1, g)=0$ for any $g \in \mathcal F$ (see e.g., \cite[Theorem 5.6.1]{O13}). We focus only on the case of RCD$^*(K,N)$ with infinite mass $m(X)=\infty$.

Let $(X,d,m, \x)$ be an RCD$^*(K,N)$ space. Note that $(X,d,m)$ becomes locally compact because of the RCD$^*(K,N)$ condition. Recall that $\partial$ denotes a cemetery point jointed to $X$ as one-point compactfication.  
Let $\mathcal A:=\{\rho \in \mathcal F_{loc}\cap C(X): \lim_{x \to \partial} \rho(x)=\infty, \{x \in X: \rho (x) \le r\}\text{ is compact for any $r>0$}\}$.
Let $B_r^\rho=\{x \in X: \rho(x) \le r\}$ and $M^\rho(r):={\rm ess}\text{-}{\rm sup}_{x \in B_r^\rho}\langle \tilde{A}\nabla \rho, \nabla \rho \rangle(x)$.
\begin{prop} \label{thm: NS2}
Let $(X,d,m,\x)$ be an RCD$^*(K,N)$ space.
Suppose Assumption \ref{asmp: NS1}\deleted{ with $\alpha_0=0$} and
\begin{align} \label{eq: DIVZ}
{\rm div} \b_i=c, \quad \text{for}\  i=1,2.
\end{align}
Assume that there exists $\rho \in \mathcal A$ so that, for any $R>0$, 
\begin{align} \label{ass: conserv1}
\lim_{r \to \infty}m(B^\rho_{R+r}){\rm Erfc}(\frac{r}{\sqrt{M^\rho(R+r)T}}) =0,
\end{align}
whereby ${\rm Erfc}(x):=\frac{2}{\sqrt{\pi}}\int_x^\infty e^{-y^2}dy$ and there exists a constant $c>0$ so that with the above $\rho$, it holds that 
\begin{align} \label{ass: conserv2}
|\b_1-\b_2||\nabla \rho|\1_{B_r^\rho} \le c(1+r) \quad \text{$m$-a.e., for any $r>0$}.
\end{align} 
Then the form $(\mathcal E, \F)$ and the dual form $(\hat{\mathcal E}, \hat{\F})$ are conservative.
\end{prop}

{\it Proof of Propostion \ref{thm: NS2}}. 
The idea of the proof is similar to the case of Euclidean diffusions discussed in \cite[Section 4]{TT11}.
We only prove the statement for $(\mathcal E, \mathcal F)$ since the dual case can be shown in the same proof. 
Let us write $m_{R}=m\1_{B^\rho_R}.$
Let $(\{S_t\}_{t \ge 0}, \mathbb P_r^x)$ denote a part process on $B_r(\x)$ of $(\{S_t\}_{t \ge 0}, \mathbb P^x)$ with $x \in B^\rho_r$\added{, which is a stopped process when it hits the boundary $\partial B_r(\x)$}. Let $T>0$. If $S_0 \in B^\rho_R$, then, by the locality of $(\E, \F)$,  
\begin{align*}
&E_r=\{\sup_{t \in [0,T]}(\rho(S_t)-\rho(S_0)) \ge r\}=\{\sup_{t \in [0,T], t<\tau_{R+r}}(\rho(S_t)-\rho(S_0)) \ge r\} \quad \text{under}\ S_0 \in B^\rho_R.
\end{align*}
Here $\tau_{R+r}=\inf\{ t \ge 0: S_t \in \partial B^\rho_{R+r}\}$. Thus, we have $\mathbb P^{m_R}(E_r)=\mathbb P_{R+r}^{m_R}(E_r).$
By \eqref{eq: DIVZ}, the form $\E$ has no killing term, which implies that the corresponding process has no inside killing.  
Therefore,
\begin{align*}
\mathbb P^{m_R}(\sup_{t \in [0,T]}(\rho(S_t)-\rho(S_0))=\infty)&=\lim_{r \to \infty}\mathbb P^{m_R}(E_r)
\\
&=\lim_{r \to \infty}\mathbb P^{m_R}_{R+r}(E_r)
\\
&=\lim_{r \to \infty}\mathbb P_{R+r}^{m_{R+r}}(\sup_{t \in [0,T], t<\tau_{R+r}}(\rho(S_t)-\rho(S_0)) \ge r ).
\end{align*}
The goal for the proof is to show 
$$\lim_{r \to \infty}\mathbb P_{R+r}^{m_{R+r}}(\sup_{t \in [0,T], t<\tau_{R+r}}(\rho(S_t)-\rho(S_0)) \ge r )=0.$$
It is easy to check that the function $\rho_{R+r}:=((\rho-(R+r)) \wedge 0)+R+r$ belongs to $\mathcal F_{B^\rho_{R+r}}:=\{f \in L^2(X,m): f|_{B^\rho_{R+r}} \in \mathcal F\}$ and 
$\rho_{R+r}=\rho$ on $B^\rho_{R+r}$.
Thus, by \eqref{eq: tight2}, we have
\begin{align*}
&\mathbb P_{R+r}^{m_{R+r}}\Bigl(\sup_{t \in [0,T], t<\tau_{R+r}}(\rho(S_t)-\rho(S_0)) \ge r \Bigr)
\\
&\le \mathbb P_{R+r}^{m_{R+r}}\Bigl(\sup_{t \in [0,T], t<\tau_{R+r}}\frac{1}{2}M_t^{[\rho]} \ge \frac{r}{4}\Bigr)+ \mathbb P_{R+r}^{m_{R+r}}\Bigl(\sup_{t \in [0,T], t<\tau_{R+r}}-\frac{1}{2}\hat{M}_t^{[\rho]}(r_T) \ge \frac{r}{4}\Bigr)
\\
&\quad +\mathbb P_{R+r}^{m_{R+r}}\Bigl(\sup_{t \in [0,T], t<\tau_{R+r}}\frac{1}{2}\hat{M}_{T-r}^{[\rho]}(r_T) \ge \frac{r}{4}\Bigr)+\mathbb P_{R+r}^{m_{R+r}}\Bigl(\sup_{t \in [0,T], t<\tau_{R+r}}\frac{1}{2}(N_t^{[\rho]}-\hat{N}_t^{[\rho]}) \ge \frac{r}{4}\Bigr).
\\
&\le \mathbb P_{R+r}^{m_{R+r}}\Bigl(\sup_{t \in [0,T], t<\tau_{R+r}}\frac{1}{2}M_t^{[\rho]} \ge \frac{r}{4}\Bigr)+ \hat {\mathbb P}_{R+r}^{m_{R+r}}\Bigl(\sup_{t \in [0,T], t<\tau_{R+r}}-\frac{1}{2}\hat{M}_t^{[\rho]} \ge \frac{r}{4}\Bigr)
\\
&\quad +\hat{\mathbb P}_{R+r}^{m_{R+r}}\Bigl(\sup_{t \in [0,T], t<\tau_{R+r}}\frac{1}{2}\hat{M}_{T-r}^{[\rho]}\ge \frac{r}{4}\Bigr)+\mathbb P_{R+r}^{m_{R+r}}\Bigl(\sup_{t \in [0,T], t<\tau_{R+r}}\frac{1}{2}(N_t^{[\rho]}-\hat{N}_t^{[\rho]}) \ge \frac{r}{4}\Bigr).
\end{align*} 
By using \eqref{ass: conserv1}, the martingale parts go to zero as $r \to \infty$ in a similar way to \cite[\S 5.7]{FOT11}, so we omit the proof. 
 We just need to estimate the zero-energy parts. 
By \replaced{the equality \eqref{eq: LZR}}{Proposition \ref{prop: NSLZ1}} and \eqref{eq: DIVZ}, we have 
$$\frac{1}{2}(N_t^{[\rho]}-\hat{N}_t^{[\rho]})=\int_0^t (\b_1-\b_2)(\rho)(S_s)ds.$$
Therefore, by \eqref{ass: conserv2}, we have
\begin{align*}
&\mathbb P_{R+r}^{m_{R+r}}\Bigl(\sup_{t \in [0,T], t<\tau_{R+r}}\frac{1}{2}(N_t^{[\rho]}-\hat{N}_t^{[\rho]}) \ge \frac{r}{4}\Bigr)
\\
&=\mathbb P_{R+r}^{m_{R+r}}\Bigl(\sup_{t \in [0,T], t<\tau_{R+r}}\int_0^t (\b_1-\b_2)(\rho)(S_s)ds \ge \frac{r}{4}\Bigr)
\\
&=\mathbb P_{R+r}^{m_{R+r}}\Bigl(|\b_1-\b_2||\nabla \rho|T \ge \frac{r}{4}\Bigr)
\\
&\le \mathbb P^{R+r}_{m_{R+r}}\Bigl(c(1+r+R)T \ge \frac{r}{4}\Bigr)
\\
& \to 0 \quad \text{as} \ r\to \infty \quad \text{if $T<\frac{1}{4c}$}.
\end{align*}
Thus, the desired result is true for $T<\frac{1}{4c}$. By using the Markov property, we can extend the result for any $T\ge 0$. Thus, we have finished the proof.
\qed

{
As a corollary, we obtain that diffusion processes constructed in Proposition \ref{thm: NS1} are conservative if ${\rm div} \b_i=c$ for $i=1,2$.
\begin{cor} \normalfont \label{ex: conserv12}
Let $X$ be an RCD$^*(K,N)$ space and $(\E, \F)$ be the Dirichlet form associated with \eqref{form: non-sym}. Under Assumption \ref{asmp: NS1} and ${\rm div} \b_i=c$ for $i=1,2$, the Dirichlet form $(\E, \F)$ is conservative. 
\end{cor}
\proof
By \eqref{asmp: basic1}, we have the following volume growth estimate: there exist positive constants $c_1, c_2$ depending only on $K$ so that $m(B_r(\x)) \le c_1e^{c_2 r^2}$ 
Take $\rho=d(\cdot, \x)$. Then we obtain
\begin{align} \label{volume: GR1}
m(B^\rho_{R+r}){\rm Erfc}(\frac{r}{\sqrt{M^\rho(R+r)T}}) \le c_3\exp\{{c_2}(R+r)^2\} \frac{\sqrt{M^\rho(R+r)T}}{r}\exp\{-\frac{r^2}{2M^\rho(R+r)T}\}.
\end{align}
Here $c_3$ is a positive constant depending only on $K$.
By $|\tilde{A}| \in L^\infty(X;m)$, it holds that $M^\rho(\cdot) \in L^\infty(X;m)$.
Therefore, R.H.S.\ of \eqref{volume: GR1} goes to zero as $r \to \infty$ and we obtain \eqref{ass: conserv1}.
 The inequality \eqref{ass: conserv2} holds immediately because $|\b_1-\b_2| \in L^\infty(X, m)$.
 \qed
 }
 \section{Examples} \label{sec: exa} 
In this section, several specific examples satisfying Assumption \ref{asmp: NS3} are given. 
There are various concrete examples of non-smooth metric measure spaces satisfying RCD conditions. See Ricci limit spaces (Sturm \cite{Sturm06, Sturm06-2}, Lott--Villani \cite{LV09}, Cheeger--Colding \cite[Example 8]{CC97}), Alexandrov spaces (Petrunin, Zhang--Zhu \cite{Pet11, ZZ10}), warped products and cones (Ketterer \cite{Ket14, Ket14a}), quotient spaces (Galaz-Garc\'ia--Kell--Mondino--Sosa \cite{GKMS17}) and infinite-dimensional spaces such as Hilbert spaces with log-concave measures (Ambrosio--Savar\'e--Zambotti \cite{ASZ09}).
Also, in \cite[Section 4]{S17}, the author explained various examples relating to the weak convergence of Brownian motions such as 
weighted Riemannian manifolds whose weighted Ricci curvature is bounded below, its pmG limit spaces, Alexandrov spaces, and Hilbert spaces with log-concave probability measures. Those examples are also available for this paper, and we omit the descriptions of those examples and we refer the reader to those references. What we discuss in this section is how to construct concrete examples of coefficients $A_n, \b^n, c_n$ satisfying Assumption \ref{asmp: NS3}. 

For any $f \in W^{1,2}(X,d,m)$, recall that we set a gradient derivation operator $\b_f$ so that
$$\b_f(g):=\langle \nabla f, \nabla g\rangle, \quad g \in W^{1,2}(X,d,m).$$
Then  we can check that $\b_f$ is an $L^2$-derivation and $|\b_f|=|\nabla f|$ (see e.g., Gigli \cite{G16} for detail).

Recall a sufficient condition for the $L^2$-strong convergence for gradient derivations according to \cite[Theorem 6.4]{AST16}. 
\begin{thm}$(${\rm \cite[Theorem 6.4]{AST16}}$)$ \label{defn: CONVD} \normalfont
Let $(X_n, d_n, m_n, \x_n)$ be a sequence of p.m.m. spaces with RCD$(K,\infty)$ condition. Assume that $(X_n,d_n,m_n, \x_n)$ converges to $(X_\infty, d_\infty, m_\infty, \x_\infty)$ in the pmG sense.  If $f_n \in W^{1,2}(m_n)$ converges strongly in $W^{1,2}$ to $f_\infty \in W^{1,2}(m_\infty)$, then $\b_{f_n}$  converges strongly in $L^2$ to $\b_{f_\infty}$.
\end{thm}
Using Theorem \ref{defn: CONVD}, we give an example of gradient derivations which is given by the resolvent of the Cheeger energy and  satisfies Assumption \ref{asmp: NS3}, according to \cite[Example 6.6]{AST16}. 
\begin{exa} {\bf (Derivation associated with resolvents)} \normalfont \label{thm: EXNS}
Let $(X_n, d_n, m_n, \x_n)$ be a sequence of p.m.m. spaces with RCD$(K,\infty)$ condition with $m_n(X_n)=1$ or RCD$^*(K,N)$. Assume that $(X_n,d_n,m_n, \x_n)$ converges to $(X_\infty, d_\infty, m_\infty, \x_\infty)$ in the pmG sense. Let $\{G^n_\lambda\}_{\lambda \ge 0}$ and $\{H_t^n\}_{t \ge 0}$ 
be the resolvent and the semigroup associated with Cheeger energy ${\sf Ch}_n.$ Let $h \in H_{\Q_+}\mathscr A_{bs}(X_\infty)$ with $h \ge 0$\deleted{, whereby $\tilde{\ \ }$ means the McShane extension of a function on $X_\infty$ to the whole space $X$}. Let $g^i_n \in L^\infty(m_n) \cap L^2(m_n)$ ($i=1,2$) 
satisfying $\sup_{n \in \N}\|g^i_n\|_\infty<\infty$ and $g^i_n$ converges to $g^i_\infty \in L^\infty(m_\infty)\cap L^2(m_\infty)$ strongly in $L^2$ for $i=1,2$. Set $A_n:=\tilde{h}|_{X_n}+a$ with a constant $a>0$, {thus 
$$\langle A_n \nabla f, \nabla f\rangle=(\tilde{h}|_{X_n}+a)|\nabla f|^2,$$
 whereby $\tilde{\ \ }$ means the McShane extension of a function on $X_\infty$ to the whole space $X$ and $|_{X_n}$ denotes the restriction of functions to $X_n$.} Let $f^i_n:=G_\lambda g^i_n$, $\b_i^n:=\b_{f^i_n}$ for $i=1,2$ and take $c_n \in L^\infty(m_n) \cap L^2(m_n)$ ($n \in \EN$) converging in $L^2$ strongly to $c_\infty$ and $c_n \ge \max\{{\rm div}\b_{f_n^1}, {\rm div}\b_{f_n^2}\}$ for all $n \in \EN$.
Then Assumption \ref{asmp: NS3} is satisfied.
\end{exa}
{\it Proof of Example \ref{thm: EXNS}}. 
By \cite[Theorem 5.7]{AH16}, it is easy to see that $A_n$ converges to $A_\infty$ in the sense of Definition \ref{defn: CONVAN}.
We only discuss $g_n^1$ and $f_n^1$, and write $g_n^1=g_n$ and $f_n^1=f_n$ for simplicity of notation in this paragraph.
Since $g_n \in L^\infty(m_n)\cap L^2(m_n)$ and $f_n:=G^n_\lambda g_n$,
we have that $G^n_\lambda g_n \in \mathcal D(\Delta_n) \subset W^{1,2}(X,d,m)\cap L^\infty(m_n)$ where $\mathcal D(\Delta_n)$ denotes the domain of the infinitesimal generator associated with ${\sf Ch}_n$. Moreover $G^n_\lambda g_n \in {\rm Lip}_{b}(X)$. In fact, since we know that, by \cite[Theorem 6.5]{AGS14b},
$$H^n_tg_n \in L^\infty(X_n,m_n), \quad {\rm Lip}(H^n_tg_n) \le \frac{||g_n||_\infty}{\sqrt{2I_{2K}(t)}},$$
where $I_K(t):=\int_0^t e^{Ks}ds.$
Noting $G_\lambda^ng_n=\int_0^\infty e^{-\lambda t}H_t^ng_n,$ we have that $\|G_\lambda^n g_n\|_\infty<\|g_n\|_\infty/\lambda$. 
and 
\begin{align*}
|G_\lambda^ng_n(x)-G_\lambda^ng_n(y)| 
&\le \int_0^\infty e^{-\lambda t}|H_t^ng_n(x)-H_t^ng_n(y)|dt
\\
&\le \int_0^\infty e^{-\lambda t}\frac{||g_n||_\infty}{\sqrt{2I_{2K}(t)}}d(x,y)dt
\\
&\le \int_0^\infty e^{-\lambda t}\frac{||g_n||_\infty}{\sqrt{\int_0^t e^{2Ks}ds}}d(x,y)dt
\\
&\le ||g_n||_\infty d(x,y)\int_0^\infty e^{-\lambda t}\frac{1}{\sqrt{\frac{e^{2Kt}-1}{2K}}}dt. 
\end{align*}
Therefore, $G^n_\lambda g_n \in {\rm Lip}_b(X_n)$ and its Lipschitz constant is uniformly bounded in $n$.
Thus, $\sup_{n \in N}\||\b_{f_n}|\|_\infty=\sup_{n \in \N}\||\nabla f_n|\|_\infty=\sup_{n \in \N}\||\nabla G_\lambda g_n|_n\|_\infty<\infty$. Since $\sup_{n \in \N}\|g_n\|_\infty<\infty$ and $g_n \to g_\infty$ in $L^2$ strongly, then $\sup_{n \in \N}\|f_n\|_\infty<\infty$, and by the result of the Mosco convergence of ${\sf Ch}_n$ to ${\sf Ch}_\infty$ \cite[Theorem  6.8]{GMS13}, we have that $f_n \to f_\infty$ strongly in $W^{1,2}$ (see also, \cite[Corollary 6.10]{GMS13}).
Therefore, by Theorem \ref{defn: CONVD}, the gradient derivation $\b_{f_n}$ converges to $\b_{f_\infty}$ strongly in $L^2$, and $\sup_{n \in \N}\||\b_{f_n}|\|_\infty<\infty$.
Since ${\rm div}\b_{f_n}=\Delta_nf_n=\lambda f_n-g_n$, we have $\sup_{n \in \N}\|{\rm div}\b_{f_n}\|_\infty<\infty$ and ${\rm div} \b_{f_n} \to {\rm div}\b_{f_\infty}$ strongly in $L^2$. 
\qed

We give another example, which is given in terms of eigenfunctions of Laplacian according to \cite[Example 6.7]{AST16}.
\begin{exa} {\bf (Derivation associated with eigenfunctions of Laplacian)} \normalfont \label{exa: EXNS}\\
Let $\mathcal X_n=(X_n,d_n,m_n,\x_n)$ be an RCD$(K,\infty)$ space for all $n \in \N$ converging to $\mathcal X_\infty=(X_\infty, d_\infty, m_\infty, \x_\infty)$ in the pmG sense. Let $u_n$ be a normalized eigenfunction $\int_{X_n}u_n^2dm_n=1$ of the generator $-\Delta_n$ associated with ${\sf Ch}_n$ with $-\Delta_nu_n=\lambda u_n$ for some $\lambda \in \R_{\ge 0}$. Assuming $K>0$, or $m_n(X_n)=1$, by \cite[Proposition 6.7]{GMS13}, we have that $-\Delta_n$ has discrete spectra $\{\lambda_n^k\}_{k=1}^\infty$ (non-decreasing order) with the eigenfunctions $\{u_n^k\}_{k \ge 0}$. By \cite[Theorem 7.8]{GMS13}, we know that $\lambda_n^k$ converges to $\lambda_\infty^k$, and $u_n^k$ converges to $u_\infty^k$ strongly in $L^2$ if the limit eigenvalue is simple (if not simple, we can extract a convergence subsequence). By ${\sf Ch}_n(u_n^k)=\lambda^k_n \to \lambda_\infty^k={\sf Ch}_\infty(u_\infty^k)$, it holds that $u_n^k$ converges to $u_\infty^k$ strongly in $W^{1,2}$, which implies $\b_{u_n^k}$ converges to $\b_{u_\infty^k}$ strongly in $L^2$ by Theorem \ref{defn: CONVD}. On the other hand, the action of the heat semigroup $H^n_tu_n^k$ is also the $k$-th eigenfunction since $\Delta_n H^n_tu_n^k=H^n_t\Delta_nu_n^k=\lambda_n^k H_t^nu_n^k$. Since $u_n^k \in W^{1,2}(X_n,d_n,m_n)$, by Lipschitz regularization of $H_t^n$  for $W^{1,2}$ elements (see e.g., \cite[Theorem 4.3]{GKO13}, the proof is available also for the case without the condition of Alexandrov spaces),  if $\{H_t\}_{t \ge 0}$ is ultra-contractive and $\|H^n_t\|_{1\to \infty}$ is uniformly bounded from above in $n$, we have that $u_n^k$ has a Lipschitz representation $\tilde{u_n^k}$ and therefore we may assume that $u_n^k$ is Lipschitz continuous and its Lipschitz constant is dominated by $$e^{-Kt}\sqrt{\|H^n_t\|_{1 \to \infty}{\sf Ch}_n(u_n^k)}=e^{-Kt}\sqrt{\|H^n_t\|_{1 \to \infty}\lambda_n^k}.$$ For instance, if $\mathcal X_n$ is RCD$^*(K,N)$ with $\sup_{n \in \N}{\rm diam}(X_n)<\infty$, or with $\inf_{x \in X_n}m_n(B_r(x))>0$ for any fixed $r>0$, $\{H_t\}_{t \ge 0}$ is ultra-contractive and $\|H^n_t\|_{1\to \infty}$ is uniformly bounded from above in $n$. This is because the volume doubling property and the Poincar\'e inequality 
 imply the Gaussian heat kernel estimate whose constants only depend on the constants appearing in the doubling and Poincar\'e inequalities, which {yields the desired} uniform ultra-contractivity.
Thus, we obtain 
$$\sup_{n \in \N}\||\b_{u_n^k}|\|_\infty=\sup_{n \in \N}\||\nabla u_n^k|\|_\infty<e^{-Kt}\sup_{n \in \N}\sqrt{\|H^n_t\|_{1 \to \infty}\lambda_n^k}<\infty.$$
Also we have 
 $$\sup_{n \in \N}\||{\rm div}\b_{u_n^k}|\|_\infty=\sup_{n \in \N}\|\Delta_n u_n^k\|_\infty<\sup_{n \in \N}\lambda_n^k<\infty.$$
Take $\b^n_1=\b_{u^k_n}$, $\b^n_2=\b_{u^{k'}_n}$ and $c_n \in L^\infty(m_n) \cap L^2(m_n)$ ($n \in \EN$) converging in $L^2$ strongly to $c_\infty$ and satisfying $c_n \ge \max\{{\rm div}\b_{u_n^k}, {\rm div}\b_{u_n^{k'}}\}$ for all $n \in \EN$.
 Thus, by taking $A_n$ as in Example \ref{thm: EXNS}, we have given an example satisfying Assumption \ref{asmp: NS3}.
\end{exa}

\section*{Acknowledgment}
The author appreciates Prof.\ Masayoshi Takeda for suggesting constructive comments for Lyons-Zheng decompositions in the non-symmetric case. 
He also expresses his great appreciation to Prof.\ Shouhei Honda for valuable comments about the convergence of derivations.
He is indebted to his wife, Anna Katharina Suzuki-Klasen for her attentive proofreading.  
Finally, he expresses his great appreciation to an anonymous referee for carefully reading the manuscript to this paper, and giving a number of insightful and constructive comments. 
This work was supported by Hausdorff Center of Mathematics in Bonn. 


\end{document}